\documentclass[a4paper, 11pt]{amsart}
\usepackage{amsfonts}
\mathchardef\emptyset="001F

\theoremstyle{plain}

\newtheorem{theorem}{Theorem}[section]
\newtheorem{lemma}[theorem]{Lemma}

\newtheorem{remark}[theorem]{Remark}
\newtheorem{definition}[theorem]{Definition}

\theoremstyle{definition}
\theoremstyle{remark}
\numberwithin{equation}{section}
\newcommand{\var}{\varphi}
\newcommand{\e}{\varepsilon}
\newcommand{\Om}{\Omega}

\newcommand{\dx}{\,dx}

\newcommand{\R}{{\mathbb R}}

\newcommand{\salt}{\noalign{\vskip .2truecm}}
\newcommand{\parent}[3]{\left #1 {#3} \right #2}

\newcommand{\graffe}[1]{\parent \{ \}{#1}}

\newcommand{\N}{{\mathbb N}}

\title[A stability result for nonlinear Neumann problems under boundary
variations]{
A stability result for nonlinear Neumann problems\\
under boundary  variations}
\author[Gianni Dal Maso]{Gianni Dal Maso}
\address[Gianni Dal Maso]{S.I.S.S.A., Via Beirut 2-4, 34014, Trieste,
Italy}
\email{dalmaso@sissa.it}
\author[Fran\c{c}ois Ebobisse]{Fran\c{c}ois Ebobisse}
\address[Fran\c{c}ois Ebobisse]{S.I.S.S.A., Via Beirut 2-4, 34014,
Trieste,
Italy}
\email{ebobisse@sissa.it}
\author[Marcello Ponsiglione]{Marcello Ponsiglione}
\address[Marcello Ponsiglione]{S.I.S.S.A., Via Beirut 2-4, 34014,
Trieste, Italy}
\email{ponsigli@sissa.it}

\begin{document}

\baselineskip3.15ex
\vskip .2truecm

\begin{abstract}
\small{In this paper we study, in dimension two,
the stability of the solutions of
some nonlinear elliptic equations
with Neumann boundary conditions, under
perturbations of the domains in the Hausdorff complementary topology.
More precisely, for every   bounded open subset $\Om$ of $\R ^2$, we
consider the
problem
$$
\begin{cases} -\,{\rm
div}\,a(x,\nabla u_\Om)\,+\,b(x, u_\Om)=0 & \mbox{ in }\quad\Om,\\
    \salt
\quad a(x,\nabla u_\Om)\cdot\nu=0 & \mbox{ on
}\quad\partial\Om, \end{cases}
$$
where $a\colon \R^2\times\R^2\to\R^2$ and $b\colon \R^2\times\R\to\R$ are two
Carath\'eodory functions which
satisfy the standard monotonicity and growth conditions
of order $p$, with $1<\,p\,\leq 2$.

Let $\Om_n$ be a uniformly bounded sequence of open sets in $\R^2$,
whose complements
$\Om_n^c$ have a uniformly bounded number of connected components.
We prove that, if
$\Om_n^c\to\Om^c$ in the Hausdorff metric and $|\Om_n|\to|\Om|$,
then $u_{\Om_n}\to u_\Om$
and $\nabla u_{\Om_n}\to \nabla u_\Om$ strongly in $L^p$.
The proof
is obtained by showing the Mosco convergence of the
Sobolev spaces $W^{1,p}(\Om_n)$ to the Sobolev space $W^{1,p}(\Om)$.
}
\end{abstract}

\maketitle

\section{Introduction}\label{Intro}
In this paper we consider nonlinear elliptic equations
with Neumann boundary conditions of the form
\begin{equation}\label{eqgen}
\begin{cases} -\,{\rm
div}\,a(x,\nabla u)\,+\,b(x, u)=0 & \mbox{ in }\quad\Om,\\
    \salt
a(x,\nabla u)\cdot\nu=0 & \mbox{ on
}\quad\partial\Om, \end{cases}
\end{equation}
where $\Om$ is a bounded open subset of $\R ^2$ and
$a\colon \R^2\times\R^2\to\R^2$ and $b\colon \R^2\times\R\to\R$ are two
Carath\'eodory functions  which satisfy suitable monotonicity,
coerciveness, and
growth conditions
(see (\ref{a1})-(\ref{a3}) below).
Our purpose is to study
the continuity of the mapping
$\Om\mapsto u_\Om$ which associates to every $\Om$ the corresponding
solution $u_\Om$. The notion of convergence we consider on the sets
$\Om$ is given by the Hausdorff complementary
topology, which is induced by the Hausdorff distance of the
complements of the sets $\Om$ (see Section~\ref{haus}).

Many examples (see \cite{daml}, \cite{mur}, \cite{delv}, \cite{con},
\cite{cor})
show that, if
we consider a uniformly bounded  sequence $\Om_h$
of open subsets of $\R^2$
which converges to an open set $\Om$
in the Hausdorff complementary topology, it may happen that
$u_{\Om_h}$ does not converge to
$u_{\Om}$.
Therefore, in this paper we investigate additional conditions on $\Om_h$
which imply that $u_{\Om_h}$ converges to
$u_{\Om}$ for any choice of the functions $a$ and $b$.
These stability results are useful in the proof of the existence of solutions
of some shape optimization problems.
More recently, similar stability results have been applied in \cite{dmro}
to study some models in fracture mechanics.

In the special case of the linear  problems
\begin{equation}\label{neumlin}
\begin{cases} -\,\Delta u\,+\,u=f & \mbox{ in }\quad\Om,\\
\salt \displaystyle\frac{\partial u}{\partial\nu}=0 & \mbox{ on
}\quad\partial\Om,  \end{cases}
\end{equation}
the stability with respect to $\Om$  was first studied by Chenais
\cite{chen} under the assumption
that the domains $\Om_h$ satisfy a uniform cone condition,
which allows to use extension operators with uniformly bounded
norms. This condition excludes
a large class of domains, like for instance domains with cracks, for
which there is an increasing interest in view of the applications to
fracture mechanics.

The stability of (\ref{neumlin}) in nonsmooth domains is studied in
\cite{buzo1}, \cite{buzo2}, \cite{buzo3} under variuos assumptions on
$\Om_h$. In \cite{chdo} the problem is studied under the hypothesis
that $\Om_h$ converges to $\Om$ in the
Hausdorff complementary topology, assuming also a uniform bound
for the lenghts of the
boundaries ${\mathcal
H}^1(\partial\Om_h)$ and for the number of connected components of
$\partial\Om_h$.
This result has been recently improved in \cite{bv1}, where the bound
on ${\mathcal
H}^1(\partial\Om_h)$ is replaced by the weaker assumption
of convergence of the two-dimensional measures of the domains, i.e.,
$|\Om_h|\to |\Om|$, which is also necessary for the stability
of (\ref{neumlin}).

In the present paper we study the stability of the nonlinear
problems (\ref{eqgen}) with respect to the
Hausdorff complementary topology,
assuming that $|\Om_h|\to |\Om|$ and that the
number of connected components of the complements
$\Om_h^c$ is uniformly bounded.
To obtain this result we reduce the problem to the
convergence  in the sense of Mosco of the
Sobolev spaces $W^{1,p}(\Om_h)$ to the Sobolev space $W^{1,p}(\Om)$,
where the exponent $p$ is related in the usual way to the growth condition
of the functions $a$
and $b$ (see Section~\ref{MC}).

The proof of this property for $1<p\le 2$ is obtained in two steps.
First, under the
same assumptions on $\Om_h$ and $\Om$,
we prove the continuity of the map $\Om\mapsto\nabla u_\Om$
for the solutions $u_\Om$ of
following nonlinear Neumann problems
\begin{equation}\label{eqgenzero}
\begin{cases} -\,{\rm
div}\,a(x,\nabla u)=0 & \mbox{ in }\quad\Om,\\
    \salt
a(x,\nabla u)\cdot\nu=0 & \mbox{ on
}\quad\partial\Om. \end{cases}
\end{equation}

This result is obtained by using the fact that the rotation by $\pi/2$
of the vectorfield $a(x,\nabla u)$ (extended to $0$ on the
complemenent $\Om^c$) is the gradient of a function
$v_\Om$ which is constant on each connected component of $\Om^c$
(Proposition~\ref{conj}). This function plays the role of the
conjugate of $u_\Om$ used in
\cite{bv2} and \cite{dmro}  in the linear case.
Another important ingredient in the proof
is a result on the stability of nonlinear Dirichlet problems
proved in \cite{bt}, which allows to show that, if each function
$v_{\Om_h}$ is
constant on each connected component of $\Om^c_h$, then their weak limit
is constant on each connected component of $\Om^c$
(Lemmas~\ref{moscostante} and~\ref{moscostante1}).

The second step in the proof of the convergence of the Sobolev
spaces $W^{1,p}(\Om_h)$ is the approximation of locally
constant functions in $\Om$ by functions belonging to
$W^{1,p}(\Om_h)$, which relies on a result obtained in \cite{bv1}.

The main difference with respect to the linear case studied in
\cite{bv1} and
\cite{bv2} is that for the nonlinear problems (\ref{eqgenzero})
we can not use the method of conformal mappings.

The hypothesis $p\le 2$ is used both in the first and in the second
step. In the case $p>2$ the stability result for (\ref{eqgen}) and
(\ref{eqgenzero}) is not true under our
hypotheses, as shown in Remarks~\ref{pgreat2} and \ref{Sobpgreat2}. 
The general form of the
limit problem will be studied in \cite{ebpo}.

In the last part of the paper we consider the case of unbounded open
sets and the case of mixed boundary value problems, with
a Dirichlet condition on a fixed part of the boundary.

{\small \tableofcontents}

\section{Notation and preliminaries}\label{Notprel}
Throughout the paper $p$ and $q$ are real numbers, with $1<p\le 2\le
q<+\infty$ and $p^{-1}+q^{-1}=1$. The scalar product of two vectors
$\xi,\,\zeta\in\R^2$ is denoted by $\xi\cdot\zeta$, and the norm of
$\xi$ by $|\xi|$.  For any $E,\,F\subset\R^2$,
$E\triangle F:=(E\setminus F)\cup(F\setminus E)$ is
the symmetric difference of $E$ and $F$, and $|E|$ is the
Lebesgue (outer) measure of~$E$.

\subsection{Deny-Lions spaces}
Given an open subset $\Om$ of $\R^2$, the {\it Deny-Lions space\/}
is defined by
$$L^{1,p}(\Om):=\{u\in L^p_{\rm loc}(\Om):\, \nabla u\in L^p(\Om,\R^2)\}.$$
It is well-known that $L^{1,p}(\Om)=W^{1,p}(\Om)$
whenever $\Om$ is bounded
and has a Lipschitz boundary.
It is also known that the set
$\{\nabla u:\, u\in L^{1,p}(\Om)\}$ is a closed subspace of
$L^p(\Om,\R^2)$.
For further properties of the spaces $L^{1,p}(\Om)$
the reader is referred to \cite{deli} and~\cite{maz1}.

In many problems it is useful to consider the following
equivalence relation in $L^{1,p}(\Om)$:
\begin{equation}\label{defreleq}
v_1\sim v_2\qquad\mbox{if and only if}
\qquad\nabla v_1=\nabla v_2\,\mbox{ a.e.\ in }\Om.
\end{equation}
The corresponding quotient space is
denoted by $L^{1,p}(\Om)/_\sim$.

\subsection{Capacity}\label{cap}
Let $1<r<\infty$. For every subset $E$ of $\R^2$, the {\it
$(1,r)$-capacity\/} of $E$ in
$\R^2$, denoted by $C_r(E)$, is defined as the infimum of
$ \int_{\R^2}(|\nabla u|^r+|u|^r)\,dx$ over the set of all functions
$u\in W^{1,r}(\R^2)$ such that $u\geq 1$ a.e.\ in a neighborhood of
$E$. If $r>2$, then $C_r(E)>0$ for every nonempty set $E$. On the
contrary, if $r=2$ there are nonempty sets $E$ with $C_r(E)=0$ (for
instance, $C_r(\{x\})=0$ for every $x\in\R^2$).

We say that a property $\mathcal P(x)$ holds $C_r$-{\it quasi everywhere}
(abbreviated as {\it $C_r$-q.e.}) in a set $E$ if it holds for all $x\in E$
except
a subset $N$ of $E$ with $C_r(N)=0$. We recall that the expression
{\it almost everywhere} (abbreviated as {\it a.e.}) refers, as usual, to
the Lebesgue measure.

A function $u\colon  E\to\overline\R$ is said to be {\it
quasi-continuous} if for every $\e$ there exists $A_\e\subset\R^2$,
with  $C_r(A_\e)<\e$,
such that the restriction of $u$ to $E\setminus A_\e$ is
continuous. If $r>2$ every quasi-continuous function is continuous,
while for $r=2$ there are quasi-continuous functions that are not
continuous.
It is well known that, for any open subset $\Om$ of
$\R^2$, every function  $u\in W^{1,r}(\Om)$ has a {\it quasi-continuous
representative}
    $\overline u\colon \Om\to\R$, which satisfies
$$
\lim_{\rho\to 0^+}\frac{1}{|B(x,\rho)|}\int_{B(x,\rho)}|u(y)-\overline
u(x)|\,dy=0\quad\mbox{for q.e.\ }x\in \Om,
$$
where $B(x,\rho)$ is the open ball with centre $x$ and radius $\rho$.
We
recall that if $u_h$ converges strongly to $u$ in $W^{1,r}(\Om)$, then
a subsequence of $\overline  u_h$ converges to $\overline  u$
pointwise $C_r$-q.e.\
on $\Om$. For these and other properties on quasi-continuous representatives
     the reader is referred to
\cite{EG}, \cite{hekima}, \cite{maz1}, \cite{ziem}.

To simplify the notation, we always identify each function $u\in
W^{1,r}(\Om)$ with its quasi-continuous representative~$\overline  u$.

\subsection{Convergence of sets}\label{haus}
We recall here  the notion of convergence in the sense of {\it Kuratowski}.
    We say that a sequence  $(C_h)$ of closed subsets of $\R^2$
converges to a closed set $C$ in the sense of Kuratowski if the
following two properties hold:

\begin{itemize}
\item[($K_1$)] for every $x\in C,$ there exists a sequence $x_h \in
C_h$ such that $x_h\to x$;
\item[($K_2$)] if $(h_k)$ is a sequence of indices converging to
$\infty$, $(x_k)$ is a sequence such that
$x_k\in C_{h_k}$ for every $k$, and
$x_k$ converges to some $x\in \R^2$, then $x\in C$.
\end{itemize}

Let us recall also that the
{\it Hausdorff distance} between two nonempty
closed subsets $C_1$ and $C_2$  of $\R^2$ is defined by
$$
d_{H}(C_1,C_2):=\max\graffe{\displaystyle\sup_{x\in C_1}{\rm
dist}\bigl(x,C_2\bigr)\,,\,\displaystyle\sup_{x\in C_2}{\rm
dist}\bigl(x,C_1\bigr)}.$$
We say that a sequence $(C_h)$ of nonempty closed subsets of $\R^2$
converges to a nonempty closed subset $C$ in the
{\it Hausdorff metric\/} if
$d_H(C_h\,,C)$ converges to $0$.

A sequence of subsets of $\R^2$  is said to be
{\it uniformly bounded\/} if there exists a bounded subset of $\R^2$ which
contains all sets of the sequence.

The convergence in the Hausdorff metric implies
the convergence in the sense of Kuratowski, while
in general the converse is false.
However, if $(C_h)$ is a uniformly bounded sequence of nonempty
closed sets in $\R^2,$
then  $(C_h)$ converges to a closed set $C$ in the Hausdorff metric
if and only if
$(C_h)$ converges to $C$ in the sense of Kuratowski.

We say that a sequence $(\Om_h)$ of open subsets of
$\R^2$ converges to an open set $\Om$ in the
{\it Hausdorff complementary topology\/}, if $d_{H}(\Om_h^c\,,\Om^c)$
converges to $0$, where $\Om_h^c$ and $\Om^c$ are the complements
of $\Om_h$ and $\Om$ in $\R^2$,
It is well-known (see, e.g., \cite[Blaschke's Selection
Theorem]{falc}) that every uniformly bounded sequence
of nonempty closed sets is compact with respect to the
Hausdorff convergence. This implies
that every uniformly bounded sequence of open sets
is compact with respect to the Hausdorff complementary topology.

Moreover, a uniformly bounded sequence of open
sets $(\Om_h)$ converges to an open set $\Om$ in the
Hausdorff complementary topology, if and only if the sequence  $(\Om_h^c)$
converges to $\Om^c$ in the sense of Kuratowski.
\subsection{The Neumann problems}
Let  $a\colon \R^2\times\R^2\to\R^2$ and $b\colon \R^2\times\R\to\R$ be two
Carath\'eodory functions that satisfy the following
assumptions: there exist $0<c_1\le c_2$, $\alpha\in L^q(\R^2)$,
and $\beta\in L^1(\R^2)$
such that, for almost every $x\in\R^2$ and for every
$\xi,\,\xi_1,\xi_2\in\R^2$
with $\xi_1\neq\xi_2$
\begin{eqnarray}
&(a(x,\xi_1)-a(x,\xi_2))\cdot(\xi_1-\xi_2)>0;&\label{a1}\\
\salt
&|a(x,\xi)|\leq \alpha(x)+c_2|\xi|^{p-1};&\label{a2}\\
\salt
&a(x,\xi)\cdot\xi\geq -\beta(x) +c_1|\xi|^p.&\label{a3}
\end{eqnarray}
We assume that $b$ satisfies the same inequalities
for every $\xi,\,\xi_1,\xi_2\in\R$.

For every open set $\Om\subset \R^2$, we consider the following
nonlinear Neumann boundary value problems, where $\nu$ denotes the
outward unit normal to $\partial\Om$:
\begin{equation}\label{eq1}
\begin{cases} -\,{\rm
div}\,a(x,\nabla u)+  b(x,u)=0 & \mbox{ in }\quad\Om,\\
    \salt
\mbox{ }a(x,\nabla u)\cdot\nu=0 & \mbox{ on
}\quad\partial\Om,\end{cases}
\end{equation}
and
\begin{equation}\label{eqaux1}
\begin{cases} -\,{\rm
div}\,a(x,\nabla v)=0 & \mbox{ in }\quad\Om\,,\\
\salt
\mbox{ }a(x,\nabla v)\cdot\nu =0 & \mbox{ on }\quad\partial\Om \,.\end{cases}
\end{equation}
A function $u$ is a solution of
(\ref{eq1}) if
\begin{equation}\label{eq2} \begin{cases} u\in W^{1,p}(\Om),\\
\salt
\displaystyle\int_{\Om}\bigl[a(x,\nabla u)\cdot\nabla z
+b(x,u)z\bigr]dx=0 & \forall z\in W^{1,p}(\Om).
\end{cases}
\end{equation}
while $v$ is a solution of  (\ref{eqaux1}) if
\begin{equation}\label{eqaux2}
\begin{cases} v \in L^{1,p}(\Om), \\
\salt
\displaystyle\int_{\Om}a(x,\nabla v)\cdot\nabla z\,\dx=0 &
\forall z\in L^{1,p}(\Om).
\end{cases}
\end{equation}
By well-known existence results for nonlinear elliptic equations
with strictly monotone operators (see, e.g., Lions \cite{lions}),
one can easily see that (\ref{eq2}) has a unique solution in
$W^{1,p}(\Om)$. Similarly one can prove that (\ref{eqaux2}) has a
solution, and that if $v_1$ and $v_2$ are solutions of (\ref{eqaux2}),
then $\nabla v_1=\nabla v_2$ a.e.\ in $\Om$. Note that problem (\ref{eqaux2})
    can be formulated in the quotient space $L^{1,p}(\Om)/_\sim$,
    where a uniqueness result holds.
\subsection{Stability of Neumann problems}\label{snp}
In order to study the stability of (\ref{eq1}) and (\ref{eqaux1})
with respect to the
variations of the open set $\Om$, we should be able to compare
two solutions defined in two different domains. For any subset $E$ of $\R^2$,
the characteristic function $1_E$ of $E$ is defined by $1_E(x):=1$
for
$x\in E$ and $1_E(x):=0$ for $x\in E^c$.
For every $u\in
L^{1,p}(\Om)$, the functions $u1_\Om$ and $\nabla u1_\Om$ are
the extensions of the functions $u$ and $\nabla u$ which vanish in $\Om^c$.
By means of these extensions, $W^{1,p}(\Om)$
will be identified with the closed linear subspace $X_\Om$ of $L^p(\R^2)\times
L^p(\R^2,\R^2)$ defined by
\begin{equation}\label{XOm}
X_\Om:=\{(u1_\Om,\nabla u1_\Om): u\in W^{1,p}(\Om)\},
\end{equation}
while the quotient space $L^{1,p}(\Om)/_\sim$ will be identified with
the closed linear subspace
$Y_\Om$ of $ L^p(\R^2,\R^2)$ defined by
\begin{equation}\label{YOm}
Y_\Om:=\{\nabla u1_\Om: u\in L^{1,p}(\Om)\}.
\end{equation}

Let $\Om$ be an open subset of $\R^2$ and let $(\Om_h)$ be a sequence
of open subsets of $\R^2$. Given a pair of
    Carath\'eodory functions $a\colon \R^2\times\R^2\to\R^2$ and
    $b\colon \R^2\times\R\to\R$
satisfying (\ref{a1})--(\ref{a3}),
let $u$ be the solution to
problem (\ref{eq1}) in $\Om$ and, for every $h$, let $u_h$ be the solution
to problem (\ref{eq1}) in~$\Om_h$.

\begin{definition}\label{stable}
We say that $\Om$ is stable for the Neumann problems (\ref{eq1})
along the sequence $(\Om_h)$ if for
every pair of functions $a$, $b$
satisfying (\ref{a1})--(\ref{a3}) the sequence
$(u_h1_{\Om_h})$ converges to $u1_\Om$ strongly in $L^p(\R^2)$
and the sequence $(\nabla u_h1_{\Om_h})$ converges to $\nabla u1_\Om$
strongly in
$L^p(\R^2,\R^2)$.
\end{definition}

Similarly, let $v$ be a solution  to
problem (\ref{eqaux1}) in $\Om$ and, for every $h$, let $v_h$ be a solution
to  problem (\ref{eqaux1}) in~$\Om_h$.

\begin{definition}\label{stableaux}
We say that $\Om$ is stable for the Neumann problems (\ref{eqaux1})
along the sequence $(\Om_h)$ if for
every function $a$
satisfying (\ref{a1})--(\ref{a3}) the sequence
$(\nabla v_h1_{\Om_h})$ converges to $\nabla v1_\Om$ strongly in
$L^p(\R^2,\R^2)$.
\end{definition}

\subsection{Mosco convergence}\label{MC}
We shall prove that the notion of stability introduced in the previous
definitions is equivalent to a notion of convergence for subspaces of a
Banach space introduced by Mosco in \cite{mosco}.

Let $\Om_h$ and $\Om$ be open subsets of $\R^2$, and let $X_{\Om_h}$
and $X_\Om$ be the corresponding subspaces defined by (\ref{XOm}).
We recall that $X_{\Om_h}$ converges to $X_\Om$ in the sense of Mosco
(see \cite[Definition 1.1]{mosco}) if the following two properties hold:
\begin{itemize}
\item[($M_1$)] for every $u\in W^{1,p}(\Om)$, there exists a sequence
$u_h\in W^{1,p}(\Om_h)$ such that $u_h1_{\Om_h}$ converges  to
$u1_\Om$  strongly in $L^p(\R^2)$ and $\nabla u_h1_{\Om_h}$ converges to
$\nabla u1_\Om$  strongly in $L^p(\R^2,\R^2)$;
\item[($M_2$)] if $(h_k)$ is a sequence of indices converging to
$\infty$, $(u_k)$ is a sequence such that
$u_k\in W^{1,p}(\Om_{h_k})$ for every $k$, and $u_k1_{\Om_{h_k}}$
converges weakly
in $L^p(\R^2)$ to a function $\phi$, while $\nabla u_k1_{\Om_{h_k}}$
converges weakly in
$L^p(\R^2,\R^2)$ to a function
$\psi$, then there exists $u\in W^{1,p}(\Om)$ such that $\phi=u1_\Om$
and $\psi=\nabla u1_\Om$ a.e.\ in $\R^2$.
\end{itemize}
\vskip .3truecm

Analogously, the convergence in the sense of Mosco of the spaces
$Y_{\Om_h}$ to $Y_\Om$ defined by (\ref{YOm}) is obtained
by using only the convergence of
the extensions of gradients, that is:
\begin{itemize}
\item[($M_1^\prime$)] for every $u\in L^{1,p}(\Om)$, there exists a sequence
$u_h\in L^{1,p}(\Om_h)$ such that $\nabla u_h1_{\Om_h}$
converges strongly to
$\nabla u1_\Om$ in
$L^p(\R^2,\R^2)$;
\item[($M_2^\prime$)] if $(h_k)$ is a sequence of indices converging to
$\infty$, $(u_k)$ is a sequence such that
$u_k\in L^{1,p}(\Om_{h_k})$ for every $k$, and $\nabla u_k1_{\Om_{h_k}}$
converges weakly  in
$L^p(\R^2,\R^2)$ to a function $\psi$, then there exists
$u\in L^{1,p}(\Om)$ such that $\psi=\nabla u1_\Om$ a.e.\ in $\R^2$.
\end{itemize}

\begin{theorem}\label{mosccont}
Let $\Om_h$ and $\Om$ be open subsets of $\R^2$, and let $X_{\Om_h}$
and $X_\Om$ be the corresponding subspaces defined by (\ref{XOm}).
Then $\Om$ is stable for the Neumann problems (\ref{eq1})
along the sequence $(\Om_h)$ if and only if $X_{\Om_h}$ converges to
$X_\Om$ in the sense of Mosco.
\end{theorem}

\begin{proof} Assume that $\Om$ is stable for the Neumann
problems (\ref{eq1})  along the sequence $(\Om_h)$. We want to prove
that $X_{\Om_h}$ converges to
$X_\Om$ in the sense of Mosco by using only the stability of the
solutions corresponding to functions $a$ and $b$ of the special form
\begin{equation}\label{gdm20}
a(x,\xi):=a_0(x)+a_1(\xi),\qquad\qquad b(x,t):=b_0(x)+b_1(t),
\end{equation}
with
\begin{equation}\label{gdm21}
       a_1(\xi):=|\xi|^{p-2}\xi, \qquad b_1(t):=|t|^{p-2}t,
      \qquad a_0\in L^q(\R^2,\R^2), \qquad  b_0\in L^q(\R^2).
\end{equation}

Let us prove
($M_1$). Given $u\in W^{1,p}(\Om)$, let
$a_0:=-|\nabla u|^{p-2}\nabla u1_\Om$ and
$b_0:=-|u|^{p-2}u1_\Om$. Then $u$ is the solution of (\ref{eq2}) in
$\Om$ with $a$ and $b$ given by (\ref{gdm20}). Let
$u_h\in W^{1,p}(\Om_h)$ be the
solution of (\ref{eq2}) in
$\Om_h$ with the same $a$ and $b$. By Definition \ref{stable}
the sequence
$(u_h1_{\Om_h})$ converges to $u1_\Om$ strongly in $L^p(\R^2)$
and $(\nabla u_h1_{\Om_h})$ converges to $\nabla u1_\Om$ strongly in
$L^p(\R^2,\R^2)$. This proves ($M_1$).

Let us prove ($M_2$). Let $(h_k)$ be a sequence of indices converging to
$\infty$ and let $(u_k)$ be a sequence, with
$u_k\in W^{1,p}(\Om_{h_k})$ for every $k$,  such that $(u_k1_{\Om_{h_k}})$
converges weakly
in $L^p(\R^2)$ to a function $\phi$, while $(\nabla
u_k1_{\Om_{h_k}})$
converges weakly in
$L^p(\R^2,\R^2)$ to a function
$\psi$. Let $a_0:=-a_1(\psi)=-|\psi|^{p-2}\psi$ and
$b_0:=-b_1(\phi)=-|\phi|^{p-2}\phi$,
let $a$ and $b$ be defined by (\ref{gdm20}), and let $u^*$ and $u^*_{h_k}$
be the solutions of problems (\ref{eq2}) in $\Om$ and $\Om_{h_k}$
respectively. By the stability assumption the sequence
$(u^*_{h_k}1_{\Om_{h_k}})$ converges to $u^*1_\Om$ strongly in $L^p(\R^2)$
and $(\nabla u^*_{h_k}1_{\Om_{h_k}})$ converges to $\nabla u^*1_\Om$
strongly in
$L^p(\R^2,\R^2)$. This implies that $a(x, \nabla u^*_{h_k}1_{\Om_{h_k}})$
converges to $a(x, \nabla u^*1_{\Om})$ strongly in $L^q(\R^2,\R^2)$
and $b(x, u^*_{h_k}1_{\Om_{h_k}})$
converges to $b(x, u^*1_{\Om})$ strongly in $L^q(\R^2)$.
Therefore
\begin{eqnarray}\nonumber
&\displaystyle \lim_{k\to\infty}\int_{\Om_{h_k}}
\bigl[a(x,\nabla u^*_{h_k})\cdot(\nabla u_{h_k}-\nabla u^*_{h_k})
+b(x,u^*_{h_k})(u_{h_k}- u^*_{h_k})\bigr]dx=
\\
&\displaystyle =\int_{\R^2}\bigl[a(x,\nabla u^*1_\Om)\cdot(\psi-\nabla
u^*1_\Om)
+b(x,u^*1_\Om)(\phi- u^*1_\Om)\bigr]dx.\label{gdm22}
\end{eqnarray}
By (\ref{eq2}) the left hand side of (\ref{gdm22}) is zero.
Therefore, using
(\ref{gdm20}) and (\ref{gdm21}) we obtain
$$
\int_{\R^2}\bigl[(a_1(\nabla u^*1_\Om)-a_1(\psi))\cdot(\psi-\nabla
u^*1_\Om)
+(b_1(u^*1_\Om)-b_1(\phi))(\phi- u^*1_\Om)\bigr]dx=0.
$$
Using the strict monotonicity of $a_1$ and $b_1$ we obtain that
$\psi=\nabla u^*1_\Om$ and $\phi=u^*1_\Om$ a.e.\ in~$\R^2$.

Conversely, assume now that $X_{\Om_h}$ converges to
$X_\Om$ in the sense of Mosco and let us prove the stability.
Let $a$ and $b$ be two Carath\'eodory
functions satisfying (\ref{a1})--(\ref{a3}) and let $u_h$ and $u$
be the solutions to problems (\ref{eq1}) in $\Om_h$ and $\Om$. The
weak convergence in $L^p(\R^2)\times L^p(\R^2,\R^2)$ of $(u_h1_{\Om_h},
\nabla u_h1_{\Om_h})$ to $(u1_{\Om},\nabla u1_{\Om})$ is a particular
case of \cite[Theorem A]{mosco}. For the reader's convenience, we
give here the simple proof.

Using $z:=u_h$ as test function in (\ref{eq2}) for $\Om_h$, from (\ref{a3}) we
obtain that $\|u_h\|_{W^{1,p}(\Om_h)}$ is bounded. Passing to a
subsequence, we have that $u_h1_{\Om_h}$ converges weakly
in $L^p(\R^2)$ to a function $\phi$, while $\nabla u_h1_{\Om_h}$
converges weakly in
$L^p(\R^2,\R^2)$ to a function
$\psi$. By ($M_2$) there exists
$u^*\in W^{1,p}(\Om)$ such that $\phi=u^*1_\Om$
and $\psi=\nabla u^*1_\Om$ a.e.\ in $\R^2$. By monotonicity for every
$v\in W^{1,p}(\Om)$ we have
\begin{eqnarray}
&\displaystyle \int_{\R^2} \bigl[a(x,\nabla v1_{\Om})
\cdot(\nabla v1_\Om-\nabla u_h1_{\Om_h})
+b(x,v1_{\Om})(v1_\Om-u_h1_{\Om_h}) \bigr]dx\ge
\label{gdm10}
\\
&\displaystyle \ge \int_{\R^2} \bigl[a(x,\nabla u_h1_{\Om_h})
\cdot(\nabla v1_\Om-\nabla u_h1_{\Om_h})
+b(x,u_h1_{\Om_h})(v1_\Om-u_h1_{\Om_h}) \bigr]dx.
\nonumber
\end{eqnarray}
By ($M_1$) there exists a sequence
$v_h\in W^{1,p}(\Om_h)$ such that $v_h1_{\Om_h}$ converges  to
$v1_\Om$
strongly in $L^p(\R^2)$ and $\nabla v_h1_{\Om_h}$ converges to
$\nabla v1_\Om$
strongly in $L^p(\R^2,\R^2)$. As $v_h-u_h\in W^{1,p}(\Om_h)$, by
(\ref{eq2}) we have
\begin{eqnarray}
&\displaystyle \int_{\R^2} \bigl[a(x,\nabla u_h1_{\Om_h})
\cdot(\nabla v1_\Om-\nabla u_h1_{\Om_h})
+b(x,u_h1_{\Om_h})(v1_\Om-u_h1_{\Om_h}) \bigr]dx=
\label{gdm11}
\\
&\displaystyle = \int_{\R^2} \bigl[a(x,\nabla u_h1_{\Om_h})
\cdot(\nabla v1_\Om-\nabla v_h1_{\Om_h})
+b(x,u_h1_{\Om_h})(v1_\Om-v_h1_{\Om_h}) \bigr]dx.
\nonumber
\end{eqnarray}
Since $a(x,\nabla u_h1_{\Om_h})$ is bounded in
$L^q(\R^2,\R^2)$ and $b(x,u_h1_{\Om_h})$ is bounded in
$L^q(\R^2)$, passing to the limit in (\ref{gdm10}) and (\ref{gdm11})
we obtain
\begin{eqnarray}\label{gdm12}
&\displaystyle \int_{\Om} \bigl[a(x,\nabla v)
\cdot(\nabla v-\nabla u^*)
+b(x,v)(v-u^*) \bigr]dx\ge \\
&\displaystyle \ge \lim_{h\to\infty}
\int_{\R^2} \bigl[a(x,\nabla u_h1_{\Om_h})
\cdot(\nabla v1_\Om-\nabla v_h1_{\Om_h})
+b(x,u_h1_{\Om_h})(v1_\Om-v_h1_{\Om_h}) \bigr]dx = 0.
\nonumber
\end{eqnarray}
Then we take $v=u^*\pm\e z$ in (\ref{gdm12}), with $z\in W^{1,p}(\Om)$ and
$\e>0$. Dividing by $\e$, and passing to the limit as $\e$ tends to $0$, we
obtain that $u^*$ satisfies (\ref{eq2}) in $\Om$. This proves that
$u^*=u$. Therefore $u_h1_{\Om_h}$ converges to $u1_\Om$ weakly
in $L^p(\R^2)$ and $\nabla u_h1_{\Om_h}$ converges to
$\nabla u1_\Om$ weakly in
$L^p(\R^2,\R^2)$.

Taking $v:=u$ in (\ref{gdm12})
we obtain that
\begin{eqnarray*}
&\displaystyle 
\int_{\R^2} (a(x,\nabla u_h1_{\Om_h})-a(x,\nabla u1_{\Om}))
\cdot(\nabla u_h1_{\Om_h} - \nabla u1_\Om)dx + 
\\ 
&\displaystyle
+ \int_{\R^2} (b(x,u_h1_{\Om_h})-
b(x,u1_{\Om})) (u_h1_{\Om_h}-u1_\Om)dx 
\end{eqnarray*}
tends to $0$ as $h\to \infty$.
Using the monotonicity of $a$ and
$b$ we conclude that each integral tends to~$0$.
The strong convergence of $(u_h1_{\Om_h}, \nabla u_h1_{\Om_h})$
is now a consequence of the
following lemma.
\end{proof}

\begin{lemma}\label{lemmastrong}
Let $(\psi_h)$ be a sequence in $L^p(\R^2,\R^2)$ converging
weakly in $L^p(\R^2,\R^2)$ to a function $\psi$. Assume that
\begin{equation}\label{gdm15}
\lim_{h\to\infty}
\int_{\R^2} (a(x,\psi_h)-a(x,\psi))
\cdot(\psi_h-\psi)dx = 0.
\end{equation}
Then $\psi_h$ converges to $\psi$ strongly in $L^p(\R^2,\R^2)$.
\end{lemma}

\begin{proof}
Various forms of this lemma have been used in the study of
Leray-Lions operators (see, e.g., \cite[Lemma 5]{BO-Mu-Pu}).
For the sake of completeness, we give here the
short proof of the present version.

Let $g_h:=(a(x,\psi_h)-a(x,\psi))\cdot(\psi_h-\psi)$. By monotonicity
we have $g_h\ge 0$ a.e.\ in $\R^2$, thus (\ref{gdm15}) implies that
$g_h$ converges to $0$ strongly in $L^1(\R^2)$. Passing to a
subsequence, we may assume that $g_h$ converges to $0$ a.e.\ in $\R^2$.
Using the Cauchy inequality, from (\ref{a2}) and (\ref{a3}) we obtain
for every $\e>0$
\begin{eqnarray*}
&c_1|\psi_h|^p\le \beta + a(x,\psi_h)\cdot\psi_h=
\beta + g_h + a(x,\psi_h)\cdot\psi + a(x,\psi)\cdot(\psi_h-\psi)\le
\\
&\displaystyle
\le g_h + \beta + \alpha|\psi| + |a(x,\psi)||\psi| +
c_2 \big(\frac{\e^q}{q} + \frac{\e^p}{p}\big)|\psi_h|^p
+ c_2\big( \frac{|\psi|^p}{p\e^p} + \frac{|a(x,\psi)|^q}{q\e^q}\big).
\end{eqnarray*}
Choosing $\e$ small enough, we obtain that there exist a constant
$c_3>0$ and a function
$\gamma\in L^1(\R^2)$ such that
\begin{equation}\label{gdm16}
       c_3|\psi_h|^p\le  g_h + \gamma.
\end{equation}

Let us fix a point $x\in \R^2$ where $\gamma(x)<+\infty$ and where
$g_h(x)$ tends to $0$. By (\ref{gdm16}) the sequence $\psi_h(x)$ is
bounded in $\R^2$, thus a subsequence (depending on $x$)
converges to a vector
$\xi\in \R^2$. By the definition of $g_h(x)$ and by the continuity of
$a(x,\cdot)$ we get
$(a(x,\xi)-a(x,\psi(x)))\cdot(\xi-\psi(x))=0$, which implies
$\xi=\psi(x)$ by (\ref{a1}). Therefore the whole sequence $\psi_h(x)$
converges to $\psi(x)$. Since this is true for a.e.\ $x\in \R^2$, the
strong convergence in $L^p(\R^2,\R^2)$ follows from (\ref{gdm16}) by
the dominated convergence theorem.
\end{proof}
\begin{remark}\label{MS->mis}
{\rm Let us observe that, if $(\Om_h)$  is a uniformly bounded sequence of
open subsets of $\R^2$ such that $X_{\Om_h}$ converges to $X_\Om$ in
the sense of Mosco,
then $|\Om_h\triangle\Om|$ converges to $0$.

Let $\Sigma\subset\R^2$ be a bounded closed set such that
$\Om_h\subset\Sigma$ for every $h$.
{}From property ($M_1$) it follows
that $\Om\subset\Sigma$. Indeed, if $\Om\setminus\Sigma\neq\emptyset$,
let $B\subset\Om\setminus\Sigma$ be an open ball
and let
$\varphi\in C^\infty_c(B)$ with $\varphi\neq 0$. By property ($M_1$),
there exists
    $u_h\in W^{1,p}(\Om_h)$ such that
$u_h1_{\Om_h}\to \varphi$ strongly in $L^p(\R^2)$. So,
$$0<\int_B|\varphi|^p\dx=\lim_{h\to\infty}\int_B|u_h1_{\Om_h}|^p\dx=0,$$
which is absurd.

Now, let $u:=1_\Om$. Since $u\in W^{1,p}(\Om)$,
by property ($M_1$) there exists $u_h\in W^{1,p}(\Om_h)$ such that
$u_h1_{\Om_h}\to 1_\Om$ strongly in $L^p(\R^2)$. As $|u_h1_{\Om_h}-
1_\Om|^p=1$ a.e.\ on $\Om\setminus\Om_h$, we have
\begin{equation}\label{limmis1}
\lim_{h\to\infty}|\Om\setminus\Om_h|
\,\leq\, \lim_{h\to\infty}\int_{\R^2}|u_h1_{\Om_h}- 1_\Om|^p\dx\, =\,0.
\end{equation}
On the other hand, up to a subsequence, $1_{\Om_h}$ converges weakly
in $L^p(\R^2,\R^2)$
to some $\phi$.
Hence from property ($M_2$), there
exists $v\in W^{1,p}(\Om)$ such that $\phi=v1_\Om$ a.e.\ in $\R^2$.
So we have that
$$\lim_{h\to\infty}|\Om_h\setminus\Om|=
\lim_{h\to\infty}\int_{\R^2} 1_{\Om_h}1_{\Sigma\setminus\Om}\dx
=\int_{\R^2}v1_\Om1_{\Sigma\setminus\Om}\dx=0,$$
which together with (\ref{limmis1}) gives
$|\Om_h\triangle\Om|\to 0$.

Note that, if the open sets $\Om_h$ are not uniformly bounded,
it is possible that
$X_{\Om_h}$ converges to $X_\Om$ in the sense of
Mosco while $|\Om_h\triangle\Om|$ does not converge to $0$.
Consider, for example, the sequence of open sets $\Om_h:=B(0,1)\cup
((B(0,h+h^{-1})\setminus\overline{B(0,h)})$
and $\Om:=B(0,1)$. We have that $|\Om_h\triangle\Om|=2\pi+h^{-2}\pi\to 2\pi$.

Let us verify that $X_{\Om_h}$ converges to $X_\Om$ in the sense of
Mosco. For every $u\in W^{1,p}(\Om)$, property  ($M_1$) is satisfied by
the sequence $u_h:=u1_\Om$.
    For property  ($M_2$), let $(h_k)$ be a sequence of indices converging to
$\infty$; let $(u_k)$ be a sequence, with
$u_k\in W^{1,p}(\Om_{h_k})$ for every $k$,  such that $u_k1_{\Om_{h_k}}$
converges weakly
in $L^p(\R^2)$ to a function $\phi$, while $\nabla u_k1_{\Om_{h_k}}$
converges weakly in
$L^p(\R^2,\R^2)$ to a function
$\psi$. We set $u:=\phi|_\Om$. As $u_k\rightharpoonup u$ weakly
in $L^p(\Om)$ and $\nabla u_k\rightharpoonup \psi|_\Om$ weakly
in $L^p(\Om,\R^2)$, we have $u\in W^{1,p}(\Om)$ and $\nabla
u=\psi|_\Om$ a.e.\ on $\Om$.
Now, for every ball $D\subset \Om^c$ we have $D\subset \Om^c_h$ for
$h$ large enough, hence
$$\int_D\phi\dx\,=\,
\lim_{k\to\infty}\int_Du_k1_{\Om_{h_k}}dx\,=\,0.$$
So, $\phi =0$ a.e.\ in $\Om^c$ and similarly also $\psi =0$ a.e.\ in $\Om^c$.
Hence, $\phi =u1_\Om$ and $\psi=\nabla u1_\Om$ a.e.\ in~$\R^2$.

Note that, in this case, $\Om_h$ converges to $\Om$ in
the Hausdorff complementary topology,
since $d_H(\Om_h^c,\Om^c)=h^{-1}\to0$. By adding a small strip whose
width tends to zero one can obtain an example with connected sets.
}
\end{remark}
The following theorem can be proved as Theorem~\ref{mosccont}.

\begin{theorem}\label{mosccontaux}
Let $\Om_h$ and $\Om$ be open subsets of $\R^2$, and let $Y_{\Om_h}$
and $Y_\Om$ be the corresponding subspaces defined by (\ref{YOm}).
Then $\Om$ is stable for the Neumann problems (\ref{eqaux1})
along the sequence $(\Om_h)$ if and only if $Y_{\Om_h}$ converges to
$Y_\Om$ in the sense of Mosco.
\end{theorem}

\begin{remark}\label{MDLmis}
{\rm If $(\Om_h)$  is a uniformly bounded sequence of
open subsets of $\R^2$ such that
$Y_{\Om_h}$ converges to $Y_\Om$ in the sense of Mosco,
then $|\Om_h\triangle\Om|$ converges to $0$.
Let $\Sigma\subset\R^2$ be a bounded set such that $\Om_h\subset\Sigma$
for every $h$. Arguing as in Remark \ref{MS->mis}, we get that
$\Om\subset\Sigma$.

Now, let $u(x):=\xi\cdot x$
with $\xi\in \R^2$ and $|\xi|=1$. Since $u\in L^{1,p}(\Om)$, by
property ($M_1'$)
there exists $u_h\in L^{1,p}(\Om_h)$ such that
$\nabla u_h1_{\Om_h}\to \nabla u1_\Om$
    strongly in $L^p(\R^2,\R^2)$. As $|\nabla u_h1_{\Om_h}-\nabla
    u1_\Om|^p=1$ a.e.\ on $\Om\setminus\Om_h$, we have
\begin{equation}\label{limmis2}
\lim_{h\to\infty}|\Om\setminus\Om_h|
    \,\leq\, \lim_{h\to\infty}\int_{\R^2}|\nabla u_h1_{\Om_h}-\nabla
u1_\Om|^p\dx\,=\,0.
\end{equation}
On the other hand, we consider the functions
$v_h\in L^{1,p}(\Om_h)$ defined by
$v_h(x):=\xi\cdot x$.
    Up to a subsequence,
$\nabla v_h1_{\Om_h}=\xi1_{\Om_h}$ converges weakly  in
$L^p(\R^2,\R^2)$ to some function $\psi$.
By property ($M_2'$), there exists a function $v\in L^{1,p}(\Om)$ such
that $\psi =\nabla v1_\Om$ a.e.\ in $\R^2$. So, it follows that
$$\xi\lim_{h\to\infty}|\Om_h\setminus\Om|
=\lim_{h\to\infty}\int_{\R^2}\nabla v_h 1_{\Om_h} 1_{\Sigma\setminus\Om}\dx
=\int_{\R^2}\nabla v 1_\Om 1_{\Sigma\setminus\Om}\dx\,=\,0,$$
which together with (\ref{limmis2}) gives  $|\Om\triangle\Om_h|\to0$.
}
\end{remark}
\begin{remark}\label{stabiauxmono}
    {\rm Theorems~\ref{mosccont} and \ref{mosccontaux} can
be applied in the following easy case. If $(\Om_h)$ is
increasing and $\Om$ is the union of the sequence, then it is easy to
see that
$X_{\Om_h}$ converges to
$X_\Om$ and $Y_{\Om_h}$ converges to
$Y_\Om$ in the sense of Mosco. Therefore every open set is stable for
the Neumann problems (\ref{eq1}) and (\ref{eqaux1})
along any increasing sequence converging to it.}
\end{remark}

\section{Mosco convergence of Deny-Lions spaces}\label{MDL}
In this section we  study the Mosco convergence of the subspaces
$Y_\Om$ introduced in~(\ref{YOm}) and
corresponding to the Deny-Lions spaces $L^{1,p}(\Omega)$.
By Theorem \ref{mosccontaux},
this is equivalent to the stability for the Neumann problems~(\ref{eqaux1}).

\begin{theorem}\label{moscograd}
Let $(\Om_h)$ be a uniformly bounded sequence of
open subsets of $\R^2$ that converges to an open set $\Om$ in the
Hausdorff complementary topology. Assume that $|\Om_h |$ converges to
$|\Om|$
and  that $\Om_h^c$ has a uniformly bounded number of
connected components. Then $\Om$ is stable for the Neumann
problems (\ref{eqaux1})
along the sequence $(\Om_h)$.
\end{theorem}

To prove Theorem~\ref{moscograd} we use the following lemmas.

\begin{lemma}\label{convcharac}
Let $(\Om_h)$ be a uniformly bounded sequence of  open subsets of $\R^2$
which converges to an open set $\Om$ in the Hausdorff complementary topology.
Assume that $|\Om_h | \to |\Om|$. Then $1_{\Om_h}\to 1_\Om$ in
measure, i.e., $|\Om_h\triangle\Om|\to 0$.
Moreover, if $\var_h\rightharpoonup \var$
weakly in $L^r(\R^2)$ for some $1<r<+\infty$, and
$\,\var_h=0$ a.e.\ in $\Om_h^c$, then $\var=0$ a.e.\ in $\Om^c$.
\end{lemma}
\begin{proof}
{}From the convergence of $\Om_h$ to $\Om$ in the Hausdorff
complementary topology
we have that $1_{\Om\setminus\Om_h}\to 0$ pointwise, hence
$|{\Om\setminus\Om_h}|\to 0$ by the dominated convergence theorem.
Since $|\Om_h | - |\Om| =|{\Om_h\setminus\Om}|-|{\Om\setminus\Om_h}|$,
and the left hand side tends to $0$ by hypothesis, we conclude that
$|{\Om_h\setminus\Om}|\to 0$ too.

Now, let $\psi\in L^\infty(\R^2)$. As $1_{\Om_h}\to 1_\Om$ strongly
in $L^s(\R^2)$ for
every $1\le s<+\infty$, we have
$$\int_{\Om^c}\var\psi \dx =\lim_{h\to\infty}\int_{\Om^c}\var_h \psi\dx
=\lim_{h\to\infty}\int_{\Om^c}1_{\Om_h}\var_h \psi\dx
=\int_{\Om^c}1_{\Om}\var\psi\dx=0,$$
which implies $\var=0$ a.e.\ in $\Om^c$.
\end{proof}

\begin{lemma}\label{moscostante}
Let $(\Om_h)$ be a uniformly bounded sequence of open
subsets of $\R^2$, converging  to an open set $\Om$ in the Hausdorff
complementary
topology. Assume that the sets
$\Om_h^c$ have a uniformly bounded number of connected components. Let
$(v_h)$ be a sequence in $W^{1,q}(\R^2)$ converging  weakly in
$W^{1,q}(\R^2)$ to a function $v$ and with $v_h = 0$ $C_q$-q.e.\ on
$\Om_h^c$. Then
$v =0$ $C_q$-q.e.\ on $\Om^c $.
\end{lemma}
\begin{proof} Let $\Delta_q$ be the $q$-Laplacian, defined by
$\Delta_q u:={\rm div}\,(|\nabla u|^{q-2}\nabla u)$.
Let $w_h$ and $w$ be the solutions of the problems
\begin{eqnarray}
    w_h\in W^{1,q}(\R^2),\quad\quad &\Delta_q w_h =
\Delta_q v\quad\mbox{in}\quad\Om_h,\quad\quad
&w_h=0\quad C_q\mbox{-q.e.\ in }\,\Om_h^c ,\label{gdm1}
\\
w\in W^{1,q}(\R^2),\quad\quad &\Delta_q w =
\Delta_q v\quad\mbox{in}\quad\Om,\quad\quad
&w=0\quad C_q\mbox{-q.e.\ in }\,\Om^c.\nonumber
\end{eqnarray}
Using a result on the stability of Dirichlet problems
proved by Bucur and Trebeschi
in \cite{bt} (see also \v Sver\'ak \cite{sver}
for the case $q=2$), we obtain that $w_h$ converges to $w$ strongly
in $W^{1,q}(\R^2)$. Taking $v_h-w_h$ as test function in (\ref{gdm1}),
which is possible since $v_h-w_h\in W^{1,q}_0(\Om_h)$ (see, e.g., \cite[Theorem
4.5]{hekima}),
we obtain
\begin{equation}
\langle \Delta_q w_h , v_h-w_h\rangle =
\langle \Delta_q v , v_h-w_h\rangle,
\label{gdm2}
\end{equation}
where $\langle\cdot,\cdot\rangle$ is the duality pairing between
$W^{-1,p}(\R^2)$ and $W^{1,q}(\R^2)$.
Passing to the limit in (\ref{gdm2}) we obtain
$$
\langle \Delta_q w , v-w\rangle =
\langle \Delta_q v , v-w\rangle,
$$
which implies $v=w$ by the strict monotonicity of $-\Delta_q$. Since,
by definition,
$w=0$ $C_q$-q.e.\ in $\Om^c$, we conclude that $v=0$ $C_q$-q.e.\ in $\Om^c$.
\end{proof}

\begin{lemma}\label{constconj} Let  $v\in
W^{1,q}(\R^2)$ and let $C_1$ and $C_2$ be two connected closed
subsets of
$\R^2$ with $C_1\cap C_2\neq\emptyset$.
If
$v$ is constant $C_q$-q.e.\ in $C_1$ and in $C_2$, then $v$ is
constant
$C_q$-q.e.\ in $C_1\cup C_2$. \end{lemma}

\begin{proof} For $q=2$ we the
reader is referred to
Proposition 2.5 in \cite{dmro}, while for $q>2$ the result follows from
the Sobolev embedding theorem, which yields
the continuity of $v$.
\end{proof}

\begin{lemma}\label{moscostante1}
Let $(\Om_h)$ be a
uniformly bounded sequence of open subsets of $\R^2$ which
converges to an open set $\Om$ in
the Hausdorff complementary topology,
and let $(v_h)$ be a sequence in
$W^{1,q}(\R^2)$, which converges to a function $v$
weakly in $W^{1,q}(\R^2)$.
Assume that
$\Om_h^c$ has a uniformly bounded number of connected components
and that every function $v_h$ is constant $C_q$-q.e.\ in each connected
component of $\Om_h^c$.
Then $v$ is constant $C_q$-q.e.\ in each connected
component of $\Om^c$.
\end{lemma}
\begin{proof}
Let $C_h^1,\ldots,C_h^{n_h}$ be the connected components of
$\Om_h^c$.  Passing to a subsequence, we can assume that
$n_h$ does not depend on~$h$, and that the sets $C_h^i$ converge in
the Hausdorff metric
to some connected sets $C^i$ as $h\to \infty$. Let us prove that $v$
is constant q.e.\ in
each~$C^i$.

This is trivial if $C^i$ contains only a single point. If $C^i$
has more than one point, there exists $r>0$ such that ${\rm
diam}(C_h^i)>2r$ for $h$ large enough. Let us prove that
the constant values $c_h^i$ taken by $v_h$ on $C_h^i$
are bounded uniformly with respect to~$h$. To this aim
let us consider a point $x_h\in C^i_h$. Since ${\rm
diam}(C_h^i)>2r$, we have $C^i_h\setminus B(x_h,r)\neq\emptyset$, and
by connectedness
\begin{equation}\label{polar}
C^i_h\cap\partial B(x_h,\rho)\neq\emptyset\qquad \hbox{ for every }
0<\rho<r.
\end{equation}
As $v_h=c_h^i$ $C_q$-q.e.\ on $C^i_h$, by using polar coordinates we
deduce from (\ref{polar}) the Poincar\'e inequality
$$
\int_{B(x_h,r)} |v_h-c_h^i|^q\dx\le M r^q \int_{B(x_h,r)} |\nabla
v_h|^q\dx,
$$
where the constant $M$ is independent of $h$, $i$, and $r$.
Since $v_h$ is bounded in $W^{1,q}(\R^2)$, it follows that
$c_h^i$ is bounded, and so it
converges (up to a subsequence) to some constant $c^i$.

To prove that $v = c^i$ $C_q$-q.e.\ on $C^i$, we fix two open balls
$B_1$ and $B_2$ with
$B_1\subset\subset B_2$, and a cut-off function $\varphi\in
C^\infty_c(B_2)$ with $\varphi=1$ in $B_1$. Then we have that
$\varphi(v_h - c_h^i) = 0$ $C_q$-q.e.\ on $(B_2\setminus C_h^i)^c$.
By Lemma \ref{moscostante} we get $\varphi(v - c^i) = 0$
$C_q$-q.e.\ in $(B_2\setminus C^i)^c$, hence  $v =c^i$ $C_q$-q.e on
$B_1\cap C^i$. As $B_1$
is arbitrary, we obtain
$v =c^i$ $C_q$-q.e.\ on $C^i$. If $C^i\cap C^j\neq\emptyset$, by
Lemma \ref{constconj} we have that $v$ is constant $C_q$-q.e.\ on $C^i\cup
C^j$. As $\Om^c$ is the union of the sets $C^i$, we conclude
that $v$ is constant $C_q$-q.e.\ on each connected
component of $\Om^c$.
\end{proof}

\begin{lemma}\label{conj}
Let $\Om$ be a bounded open subset of $\R^2$ and let $u$ be a
solution of problem (\ref{eqaux1}). Let $R$ be the rotation on
$\R^2$ defined by $R(y_1,y_2):=(-y_2,y_1)$.
Then there exists a unique function $v \in
W^{1,q}(\R^2)$ such that
$\nabla v=Ra(x,\nabla u)1_\Om$ a.e.\ in $\R^2$.
Moreover
$v$ is constant $C_q$-q.e.\ on each connected component of
$\Om^c$.
\end{lemma}
\begin{proof}
We consider the vector field $\Phi\in L^q(\R^2,\R^2)$
defined by $\Phi:=a(x,\nabla u)1_\Om$. By (\ref{eqaux2}) we have ${\rm
div}\,\Phi=0$ in $\mathcal D'(\R^2)$, hence
${\rm rot}(R\Phi)=0$ in
$\mathcal D'(\R^2)$. Since $\Om$ is bounded,
there exists a potential $v\in W^{1,q}(\R^2)$ such that
$\nabla v=R\Phi$ a.e.\  in $\R^2$ and $v=0$ a.e.\ in the interior of
the unbounded connected component of $\Om^c$.

Given a connected component $C$ of $\Om^c$, it remains
to prove that $v$ is constant $C_q$-q.e.\ on $C$.
For every $\e>0$ let
$C_\e:=\{x\in\R^2: {\rm dist}(x,C)<\e\}$,
and let $u_\e$ be a solution of problem (\ref{eqaux1}) in
$\Om_\e:=\Om\setminus \overline C_\e$. Let $v_\e$ be the unique function in
$W^{1,q}(\R^2)$ such that
$\nabla v_\e=Ra(x,\nabla u_\e)1_{\Om_\e}$ a.e.\ in $\R^2$.
By Remark \ref{stabiauxmono}, $\nabla
u_\e$ converges to $\nabla u$  strongly in $L^p(\R^2,\R^2)$ and so $v_\e$
converges to $v$  strongly in $W^{1,q}(\R^2)$.
By construction $\nabla v_\e=0$ in $C_\e$. As $C_\e$ is a connected
open set containing $C$, we have that $v_\e$ is constant $C_q$-q.e.\ on
$C$.  Since a subsequence of $v_\e$
converges to $v$ $C_q$-q.e.\ on $\R^2$, we conclude that
$v$ is constant $C_q$-q.e.\ on $C$.
\end{proof}

\begin{proof}[Proof of Theorem~\ref{moscograd}.]
Let $a\colon \R^2\times\R^2\to\R^2$  be a Carath\'eodory
function satisfying (\ref{a1})--(\ref{a3}) and let $u_h$ and $u$
be solutions to problems (\ref{eqaux1}) in $\Om_h$ and $\Om$.
Taking $u_h$ as test function in (\ref{eqaux2}) in $\Om_h$ and using
(\ref{a3}) we obtain that $\nabla u_h1_{\Om_h}$ is
bounded in $L^p(\R^2,\R^2)$. By (\ref{a2}) we obtain also that
$a(x,\nabla u_h)1_{\Om_h}$
is bounded in $L^q(\R^2)$. Passing to a
subsequence, we may assume that
$\nabla u_h1_{\Om_h}\rightharpoonup\Psi$
weakly in $L^p(\R^2,\R^2)$ and
$a(x,\nabla u_h)1_{\Om_h}\rightharpoonup\Phi$
weakly in $L^q(\R^2,\R^2)$. By (\ref{eqaux2}) we have ${\rm div}(
a(x,\nabla u_h)1_{\Om_h})=0$ in ${\mathcal D}'(\R^2)$, hence
${\rm div}\,\Phi=0$ in ${\mathcal D}'(\R^2)$.

If $\Om'\subset\subset\Om$, by the Hausdorff complementary convergence
we have $\Om'\subset\subset\Om_h$ for $h$ large enough.
Since the set of gradients of  functions of $L^{1,p}(\Om')$
is closed in $L^p(\Om',\R^2)$, the vectorfield $\Psi$ is the
gradient of a function of $L^{1,p}(\Om')$. As $\Om'$ is arbitrary,
we can construct
$u^* \in L^{1,p}(\Om)$ such that
$\Psi=\nabla u^*$ a.e.\ in $\Om$. On the other hand, by Lemma
\ref{convcharac} we have
$\Psi=0$ a.e.\ in $\Om^c$, hence $\Psi=\nabla u^*1_\Om$ a.e.\ in $\R^2$.

Let us prove that $\Phi=a(x,\nabla u^*)1_\Om$ a.e.\ in $\R^2$. By
Lemma \ref{convcharac}
it is enough to prove the equality in every
open ball $B\subset\subset\Om$. Note that by the Hausdorff
complementary convergence we have $B\subset\subset\Om_h$ for $h$
large enough. By adding suitable constants, we may
assume that the mean values of $u_h$ and $u^*$ on $B$ are zero. Thus
the Poincar\'e inequality and the Rellich theorem imply that $u_h$
converges to $u^*$ strongly in $L^p(B)$.

Let $z\in W^{1,p}(B)$
and let $\var\in C^\infty_c(B)$ with $\var\ge 0$. For $h$ large
enough we have $B\subset\subset\Om_h$, thus by monotonicity we
have
\begin{equation}\label{gdm30}
\int_B (a(x,\nabla z)-a(x,\nabla u_h))\cdot(\nabla z-\nabla
    u_h)\var dx\ge0.
\end{equation}
By (\ref{eqaux2}) we have also
$$
\int_B a(x,\nabla u_h)\cdot\nabla (( z-
u_h)\,\var) dx =0,
$$
which, together with (\ref{gdm30}), gives
\begin{equation}\label{gdm31}
\int_B a(x,\nabla z)\cdot\nabla (( z-
u_h)\var) dx
    -\int_B (a(x,\nabla z)-a(x,\nabla u_h))\cdot\nabla \var
    \,(z-u_h)dx\ge 0.
\end{equation}
We can pass to the limit in each term of (\ref{gdm31}) and we get
\begin{equation}\label{gdm32}
\int_B a(x,\nabla z)\cdot\nabla (( z-
u^*)\var) dx
-\int_B (a(x,\nabla z)-\Phi)\cdot\nabla \var
    \,(z-u^*)dx\ge0.
\end{equation}
As ${\rm div}\,\Phi=0$ in ${\mathcal D}'(B)$, we have
\begin{equation}\label{gdm33}
\int_B \Phi\cdot\nabla (( z-
u^*)\var) dx =0.
\end{equation}
{}{}From (\ref{gdm32}) and (\ref{gdm33}) we obtain
$$
\int_B (a(x,\nabla z)-\Phi)\cdot(\nabla z-\nabla
    u^*)\,\var dx\ge0.
$$
As $\var$ is arbitrary, we get
$(a(x,\nabla z)-\Phi)\cdot(\nabla z-\nabla
    u^*)\ge0$
    a.e.\ in $B$.
In particular, taking $z(x):= u^*(x)\pm\e \xi\cdot x$, with $\xi\in\R^2$ and
$\e>0$, we obtain
$\pm (a(x,\nabla  u^*\pm\e
\xi)-\Phi)\cdot\xi\ge0$ a.e.\ in $B$.
As $\e$ tends to zero we get
$(a(x,\nabla  u^*)-\Phi)\cdot\xi=0$ a.e.\ in $B$,
which implies that $a(x,\nabla  u^*)=\Phi$ a.e.\ in $B$ by the
arbitrariness of $\xi$.

Let us prove now that
\begin{equation}\label{eqlimit}
\int_\Om a(x,\nabla u^*)\cdot\nabla z\dx=0\qquad\forall z\in
L^{1,p}(\Om).
\end{equation}
By Lemma \ref{conj} for every $h$ there exists
$v_h\in W^{1,q}(\R^2)$
such that $\nabla v_h= Ra(x, \nabla u_h)1_{\Om_h}$ a.e.\ in $\R^2$.
Moreover $v_h$ is constant $C_q$-q.e.\ on each connected component of
$\Om_h^c$. As $a(x, \nabla u_h)1_{\Om_h}$ converges to
$a(x, \nabla u^*)1_{\Om}$ weakly in $L^q(\R^2,\R^2)$, there exists a
function $v\in W^{1,q}(\R^2)$ such that $ v_h \rightharpoonup  v$
weakly in $W^{1,q}(\R^2)$ and  $\nabla v=R
a(x, \nabla u^*)1_\Om$ a.e.\ in $\R^2$.
So, we have to prove that
\begin{equation}\nonumber
\int_\Om R\nabla v\cdot \nabla z\dx=0\qquad\forall z\in
L^{1,p}(\Om).
\end{equation}
{}From the Lemma  \ref{moscostante1} it follows that $v$ is constant
$C_q$-q.e.\ on the connected components of $ \Om^c$.
By \cite[Theorem 4.5]{hekima} we can
approximate  $v$ strongly in $W^{1,q}(\R^2)$ by  a sequence of functions
$v_h \in C^\infty_c(\R^2)$
that are constant in suitable neighborhoods $U^i_h$ of each
connected component $C^i$ of
$\Om^c$.
Let $z\in  L^{1,p}(\Om)$ and $z_h \in   W_0^{1,p}(\Om)$ such
that $ z_h=z$ in $\Om \setminus \bigcup_{i} U^i_h$.
Then, we have
\begin{equation}\label{lasteqn}
\int_\Om R\nabla v_h \cdot \nabla  z\dx = \int_\Om R \nabla v_h
\cdot \nabla z_h \dx = 0,
\end{equation}
where the last equality follows from the fact that the vector field
$R\nabla v_h$ is divergence free.
Then, passing to the limit in (\ref{lasteqn}) for $h\to\infty$, we
get
\[
    \int_\Om a(x,\nabla u^*)\cdot \nabla z\dx=
    -\int_\Om R\nabla v \cdot \nabla z
dx=0.\]
So $u^*$ is a solution of (\ref{eqaux2}) in $\Om$, hence
$\nabla
u^*=\nabla u$ a.e.\ in $\Om$ by uniqueness of the gradients.
This implies that
$\nabla u_h1_{\Om_h}$ converges to $\nabla u1_{\Om}$
weakly in $L^p(\R^2,\R^2)$ and $a(x,\nabla u_h)1_{\Om_h}$
converges to $a(x,\nabla u)1_{\Om}$ weakly in $L^q(\R^2,\R^2)$.
Since $|\Om_h \triangle \Om|$ tends to $0$ by Lemma \ref{convcharac},
from
the identity
$a(x,\nabla u_h1_{\Om_h})=a(x,\nabla u_h)1_{\Om_h}+
a(x,0)1_{\Om_h^c}$ we conclude also that
$a(x,\nabla u_h1_{\Om_h})$ converges to $a(x,\nabla u1_{\Om})$ weakly
in $L^q(\R^2,\R^2)$.

To prove the strong convergence, we consider the integral
$$
I_h:=\int_{\R^2} (a(x,\nabla u_h1_{\Om_h})-a(x,\nabla u1_{\Om}))
\cdot(\nabla u_h1_{\Om_h} - \nabla u1_\Om)dx.
$$
Since by (\ref{eqaux2})
$$
\int_{\Om_h} a(x,\nabla u_h)\cdot\nabla u_h \dx=0
\qquad\mbox{and}\qquad
\int_{\Om} a(x,\nabla u)\cdot\nabla u \dx=0,
$$
we have
$$
I_h=-\int_{\Om} a(x,\nabla u_h1_{\Om_h})
\cdot\nabla u\,dx -
\int_{\Om_h}a(x,\nabla u1_{\Om})\cdot\nabla u_h1_{\Om_h} \dx.
$$
Therefore
\begin{equation}\label{gdm40}
\lim_{h\to\infty} I_h=-2\int_{\Om} a(x,\nabla u)
\cdot\nabla u\dx
= 0,
\end{equation}
where the last equality can be deduced from (\ref{eqaux2}).
The strong convergence in $L^p(\R^2,\R^2)$ of
$\nabla u_h1_{\Om_h}$ to $\nabla u1_{\Om}$
follows now from (\ref{gdm40}) and from Lemma~\ref{lemmastrong}.
\end{proof}
\begin{remark}\label{pgreat2}
{\rm As we said in the introduction, in the case $p>2$ the stability
result for problem (\ref{eqaux1}) is not true under our hypotheses. Indeed,
let us  consider
$$S:=[1,3]\times\{0\},\quad
S_h:=([1,2-1/h]\cup [2 + 1/h,3])\times\{0\},$$
$$\Om:= B(0,3)\setminus(\overline{B(0,1)}\cup S)
\quad\mbox{ and }\quad\Om_h:= B(0,3)\setminus(\overline{B(0,1)}\cup S_h).$$
Let $\var\in C^\infty_c(0,\infty)$ be such that
$\var(\rho)=\rho^{-p}$ for $1\leq \rho\leq 3$.
We set
$$
a(x,\xi):=|\xi|^{p-2}\xi - \var(|x|)Rx,
$$
where $R$ is the rotation by $\pi/2$ defined by $Rx:=(-x_2,x_1)$.
Let $u_h$ and $u$ be solutions of problems
(\ref{eqaux1})  in $\Om$ and  $\Om_h$, with
$\int_{\Om_h} u_h\dx=\int_\Om u\dx=0$.
For every $x\in\Om$, let $0<\theta(x)<2\pi$ be the angle between $x$ and the
positive $x_1$-axis. As $\nabla \theta(x) =Rx/|x|^2$, we have
that $u=\theta -\pi$ in $\Om$.

If the open set $\Om$
  were stable for problem (\ref{eqaux1}) along the sequence $(\Om_h)$,
then $\nabla u_h$ would converge
  strongly to $\nabla u$ in $L^p(\Om,\R^2)$.
By the Poincar\'e inequality we would have
that $u_h$ converges strongly to $u$ in $W^{1,p}(\Om)$.

For every $v\in W^{1,p}(\Om)$, let
  $v^+$ and $v^-$ be the upper and lower traces of $v$ on $S$,
defined by
\begin{equation}\label{traces}
v^+(x):=\lim_{\substack{y\to x\\{y_2}>0}}v(y)\quad\mbox{ and }\quad
v^-(x):=\lim_{\substack{y\to x\\{y_2}<0}}v(y).
\end{equation}
 From the convergence of $u_h$ to $u$ in $W^{1,p}(\Om)$,
we obtain that $u_h^+\to u^+$ and $u_h^-\to u^-$ uniformly on $S$
(recall that $p>2$ here). Since $u_h^+(2,0)=u_h^-(2,0)$ by the continuity of $u_h$,
 we obtain $u^+(2,0)=u^-(2,0)$, which contradicts the
  fact that $u^+(2,0)=-\pi$ and $u^-(2,0)=\pi$, being $u={\theta -\pi}$. 
}
\end{remark}
\section{Mosco convergence of Sobolev spaces}
In this section we study
the convergence in the sense of Mosco of
the subspaces $X_{\Om}$ introduced in (\ref{XOm}) and corresponding to
the Sobolev spaces $W^{1,p}(\Om)$.
The convergence of $X_{\Om_h}$ to $X_{\Om}$
will be obtained from the  convergence of
$Y_{\Om_h}$ to $Y_{\Om}$ and from
the following  approximation theorem for
functions which are locally constant on
the limit open set $\Om$.

\begin{theorem}\label{apcostanti}
Let  $(\Om_h)$ be a uniformly
bounded sequence of  open subsets of $\R^2$ which converges to an
open set $\Om$ in the
Hausdorff complementary topology. Assume that $|\Om_h|$ converges to
$|\Om|$ and
that  $\Om_h^c$ is connected  for every $h$.
Then for every $u\in W^{1,p}(\Om)$  with $\nabla u=0$ a.e.\ in $\Om$,
there exists a sequence
$u_h\in W^{1,p}(\Om_h)$
such that $u_h1_{\Om_h}$ converges to $u1_{\Om}$ strongly in $L^p(\R^2)$ and
$\nabla u_h1_{\Om_h}$ converges to $0$ strongly in $L^p(\R^2,\R^2)$.
\end{theorem}

The proof of this theorem is  postponed.
We are now in a position to state the main result of the paper.

\begin{theorem}\label{moscosobolev}
    Let $(\Om_h)$ be a uniformly bounded sequence of
open subsets of $\R^2$ which converges to an open set $\Om$ in
the Hausdorff
complementary topology, with $|\Om_h |$ converging to $|\Om|$. Assume that
$\Om_h^c$ has a uniformly bounded number of connected components.
Then $X_{\Om_h}$ converges to
$X_\Om$ in the sense of Mosco.
\end{theorem}

To prove Theorem \ref{moscosobolev} we need the following
localization lemma.

\begin{lemma}\label{Mloc->Mglob}
Let $(\Om_h)$ be a uniformly bounded sequence of open
subsets of $\R^2$, and let $\Om$ be a bounded open subset of $\R^2$.
Assume that for every $x\in \R^2$ there exists $\e>0$ such that the
sequence $X_{B(x,\e)\cap\Om_h}$ converges to $X_{B_\e(x)\cap\Om}$
in the sense of Mosco. Then $X_{\Om_h}$ converges to $X_{\Om}$
in the sense of Mosco.
\end{lemma}
\begin{proof}
Condition ($M_2$) is easy, and condition ($M_1$) can
be obtained by using a partition of unity.
\end{proof}

\begin{proof}[Proof of Theorem \ref{moscosobolev}.] {\it Step 1.}
       Assume that $\Om^c_h$ is connected for every $h$.
Let us prove ($M_2$).  Let $(h_k)$ be a sequence of indices converging to
$\infty$, $(u_k)$ be a sequence, with
$u_k\in W^{1,p}(\Om_{h_k})$ for every $k$,  such that $u_k1_{\Om_{h_k}}$
converges weakly
in $L^p(\R^2)$ to a function $\phi$, while $\nabla u_k1_{\Om_{h_k}}$
converges weakly in
$L^p(\R^2,\R^2)$ to a function $\psi$.
{}From Lemma \ref{convcharac} it follows that $\phi$ and $\psi$
vanish a.e.\ in $\Om^c$.

Let $\Om'\subset\subset\Om$ be an open set. By the Hausdorff
complementary convergence we have $\Om'\subset\subset\Om_h$ for $h$
large enough.
So, $u_k|_{\Om'}$ converges weakly to $\phi|_{\Om'}$ in $L^p(\Om')$ and
$\nabla u_k|_{\Om'}$ converges weakly to $\psi|_{\Om'}$ in $L^p(\Om',\R^2)$.
Hence $\phi|_{\Om'}\in W^{1,p}(\Om')$ and
$\psi|_{\Om'}=\nabla\phi|_{\Om'}$ in $\Om'$.
{}From the arbitrariness of $\Om'$, it follows that the function
$u:=\phi|_\Om$ belongs to $ W^{1,p}(\Om)$,
    $\phi= u1_\Om$ and $\psi =\nabla u1_\Om$ a.e.\ in $\R^2$.

Now let us prove ($M_1$). Let $u\in W^{1,p}(\Om)$.
We write $\Om:=\bigcup_{i=1}^m\Om_i$,
where $1\leq m\leq\infty$ and $(\Om_i)$ is the family of connected
components of $\Om$.
Since the set of functions $u$ satisfying ($M_1$)
is a closed linear subspace of $W^{1,p}(\Om)$, by a density argument
it is enough to prove ($M_1$) when $u$ belongs to $L^\infty(\Om)$
and vanishes on all connected
components of $\Om$ except one. By renumbering the sequence $(\Om_i)$,
we may assume that  $u$ vanishes on $\Om_i$ for every $i\geq 2$.

{}From Theorem \ref{moscograd}
on the convergence of $Y_{\Om_h}$
to  $Y_\Om$ in the sense of Mosco, there exists
a sequence $z_h \in
L^{1,p}(\Om_h)$ such that $\nabla z_h1_{\Om_h}$
converges strongly to $\nabla u1_\Om $ in $L^p(\R^2,\R^2).$
Let us fix a nonempty open set $A_0\subset\subset\Om_1$.
We can assume that $\int_{A_0}z_hdx=\int_{A_0} udx$.
For every smooth connected open set $A$, with
$A_0\subset\subset A\subset\subset\Om_1$,
by the Poincar\'e inequality we have that $z_h|_A\to  u|_A $ strongly
in $W^{1,p}(A)$.

We consider now $w_h:=(-\|u\|_\infty)\vee z_h\wedge \|u\|_\infty$.
We have that $w_h|_A\to  u|_A $ strongly in $W^{1,p}(A)$ for every
open set
$A\subset\subset\Om_1$.
Moreover, for every open set $E\subset\subset\Om$
the function $w_h|_E$ belongs to
$W^{1,p}(E)$ for $h$ large enough. As $\|w_h\|_\infty\leq\|u\|_\infty$
and $|\nabla w_h|_E|\leq |\nabla z_h|_E|$, the sequence $(w_h|_E)$
is bounded in $W^{1,p}(E)$. By the Rellich
theorem, there exists $w\in W^{1,p}(\Om)$
such that $w_h|_E$ converges to $w|_E$
strongly in $L^p(E)$ for every open set $E\subset\subset\Om$
with smooth boundary. As
$w_h|_A\to  u|_A $ strongly in $W^{1,p}(A)$ for every open set
$A\subset\subset\Om_1$, we have that
$w=u$ a.e. in $\Om_1$.

For every open set $E\subset\subset\Om\setminus\Om_1$,
since $|\nabla w_h|\leq |\nabla z_h|$
and $\nabla z_h|_E\to\nabla u|_E=0$ strongly in $L^p(E,\R^2)$,
we have $\nabla w_h|_E\to 0$ strongly in $L^p(E,\R^2)$.
Therefore we get $\nabla w=0=\nabla u$ a.e. in
$\Om\setminus\Om_1$, which together with
the result obtained in $\Om_1$ implies that
$\nabla w_h|_E$ converges to $\nabla u|_E$
strongly in $L^p(E,\R^2)$ for every $E\subset\subset\Om$. In particular,
we obtain that the function $u-w$ is locally constant in $\Om$.

We claim  that $\nabla w_h1_{\Om_h}$ converges to $\nabla w1_\Om$ strongly in
$L^p(\R^2,\R^2)$.
Indeed for every $E\subset\subset\Om$
we have
$
\|\nabla w_h1_{\Om_h}-\nabla u1_\Om\|_{L^p(\R^2,\R^2)} \leq
\|\nabla w_h1_{E}-\nabla u1_{E}\|_{L^p(E,\R^2)}
+\, \|\nabla z_h1_{\Om_h}1_{E^c}\|_{L^p(\R^2,\R^2)}
   +\, \|\nabla u1_{\Om\setminus E}\|_{L^p(\R^2,\R^2)}$.
Hence
$$\limsup_{h\to\infty}\|\nabla w_h1_{\Om_h}-
\nabla u1_\Om\|_{L^p(\R^2,\R^2)}
\le 2\|\nabla u1_{\Om\setminus E}\|_{L^p(\R^2,\R^2)},
$$
and by letting $E\nearrow\Om$ we prove the claim. In a similar way,
we obtain also that $ w_h1_{\Om_h}$ converges to $w1_\Om$
strongly in $L^p(\R^2)$.

Since $u-w$ is locally constant in $\Om$, from Theorem \ref{apcostanti}
there exists $v_h\in W^{1,p}(\Om_h)$ such that $v_h1_{\Om_h}\to
(u-w)1_\Om$ strongly in $L^p(\R^2)$
    and $\nabla v_h1_{\Om_h}\to 0$ strongly in   $L^p(\R^2,\R^2)$.
Therefore $w_h+v_h\in W^{1,p}(\Om_h)$,
$(w_h+v_h)1_{\Om_h}\to u1_\Om$ strongly in $L^p(\R^2)$,
and $\nabla (w_h+v_h)1_{\Om_h}\to \nabla u1_\Om$ strongly in
$L^p(\R^2,\R^2)$, which give property ($M_1$).
\vskip .2truecm

\noindent {\it Step 2.} We now remove the hypothesis that
$\Om^c_h$ is connected.
Let $C_h^1,\ldots,C_h^{n_h}$ be the connected
components of
$\Om_h^c$.
Passing to a subsequence we can assume that $n_h$ does not depend on
$h$ and that the sets $C_h^i$ converge in the Hausdorff metric to some
connected sets $C^i$ as $h\to\infty$.
Let $C^{i_1},\ldots,C^{i_d}$
be those $C^i$ having at least two points. We set
$$\Om^*:=\Bigl(\bigcup_{j=1}^dC^{i_j}\Bigr)^c\quad\mbox{ and }
\quad \Om^*_h:=\Bigl(\bigcup_{j=1}^dC^{i_j}_h\Bigr)^c.$$
We have that $\Om\subset\Om^*$, $\,\Om_h\subset\Om^*_h$, and, by
construction, $\Om^*_h$ converges in the Hausdorff
complementary topology to $\Om^*$ and $|\Om^*_h|\to|\Om^*|$ (because
$|C_h^i|\to 0$ if $i\neq i_1,\ldots,i_d$).
There exists some $\eta>0$ such that ${\rm diam}\,C^{i_j}>\eta$
for every $j$, hence ${\rm diam}\,C^{i_j}_h>\eta$ for $h$ large enough.
    Let us  observe that, for every $x\in\R^2$, the sequence
$B(x,\eta /2)\cap\Om^*_h$ converges to $B(x,\eta /2)\cap\Om^*$ in the
Hausdorff complementary  topology  and also
$|B(x,\eta /2)\cap\Om^*_h|\to |B(x,\eta /2)\cap\Om^*|$ (by Lemma
\ref{convcharac}).
As ${\rm diam}\,C^{i_j}_h>\eta$, it is easy to see that
$(B(x,\eta /2)\cap\Om^*_h)^c$
is connected for $h$ large enough.
So, from  Step 1 we obtain
    the Mosco convergence of $X_{B(x,\eta /2)\cap\Om^*_h}$ to
$ X_{B(x,\eta /2)\cap\Om^*}$.
Now, using  Lemma \ref{Mloc->Mglob}  we get the Mosco convergence of
$X_{\Om^*_h}$ to $X_{\Om^*}$.

As $\Om_h^*\setminus\Om_h$ is the union of a uniformly
bounded number of sets with diameter tending to $0$,
using the fact that $1<p\leq 2$, we deduce
that $C_p(\Om_h^*\setminus\Om_h)\to 0$.
Let us show that this implies that $X_{\Om_h}$ converges to
$X_\Om$ in the sense of Mosco.

For property ($M_1$), let $u\in W^{1,p}(\Om)$. Since the set
$\Om^*\setminus\Om$ is finite, we have that $C_p(\Om^*\setminus\Om)=
0$. Hence, $u\in W^{1,p}(\Om^*)$.
So, there exists $u_h^*\in W^{1,p}(\Om_h^*)$ such that
$u_h^*1_{\Om_h^*}$ converges to $u1_{\Om^*}$ strongly
in $L^p(\R^2)$ and $\nabla u_h^*1_{\Om_h^*}$ converges to $\nabla
u1_{\Om^*}$  strongly in
$L^p(\R^2,\R^2)$. Setting $u_h=u_h^*|_{\Om_h}$, we obtain that $
u_h1_{\Om_h}$
converges to $u1_{\Om}$ strongly
in $L^p(\R^2)$ and $\nabla u_h1_{\Om_h}$ converges to $\nabla
u1_{\Om}$  strongly in
$L^p(\R^2,\R^2)$, and so property ($M_1$) holds.

Let us prove property ($M_2$). Let $(h_k)$ be a sequence of indices
converging to
$\infty$, $(u_k)$ be a sequence, with
$u_k\in W^{1,p}(\Om_{h_k})$ for every $k$,  such that $u_k1_{\Om_{h_k}}$
converges weakly
in $L^p(\R^2)$ to a function $\phi$, while $\nabla u_k1_{\Om_{h_k}}$
converges weakly in
$L^p(\R^2,\R^2)$ to a function $\psi$.
{}From Lemma \ref{convcharac} it follows that $\phi$ and $\psi$
vanish a.e.\ in $\Om^c$.

As $C_p(\Om_{h_k}^*\setminus\Om_{h_k})\to 0$, there exists a sequence
$\var _k\in W^{1,p}(\R^2)$
converging strongly to $0$ in $W^{1,p}(\R^2)$ such that $\var _k=1$ a.e.\ in
a neighborhood of $\Om^*_{h_k}\setminus\Om_{h_k}$.

We set
$$u_k^*:=u_k(1-\var_k).$$
Note that $u_k^*\in W^{1,p}(\Om^*_{h_k})$ and that
$u_k^*1_{\Om^*_{h_k}}$
converges weakly
in $L^p(\R^2)$ to $\phi$, while $\nabla u_k^*1_{\Om^*_{h_k}}$
converges weakly in
$L^p(\R^2,\R^2)$ to $\psi$. So, from the Mosco convergence of
$X_{\Om^*_h}$ to $X_{\Om^*}$, it follows that
    there exists $u^* \in W^{1,p}(\Om^*)$ such that $\phi= u^*1_{\Om^*}$
and $\psi=\nabla u^*1_{\Om^*}$
a.e.\ in $\R^2$. By setting $u=u^*|_\Om$, we get that $\phi=u1_\Om$
and $\psi=\nabla u1_\Om$
a.e.\ in $\R^2$ and the proof of ($M_2$) is complete.
\end{proof}

The rest of this section is devoted to the proof of \mbox{Theorem
\ref{apcostanti}.}
To this aim we will need some
preliminary results.

\begin{lemma}\label{compopalla2}
Let $(\Om_h)$ be a uniformly bounded sequence of open subsets
of $\R^2$, which converges to an open set $\Om$ in the Hausdorff
complementary topology, with
$|\Om_h|$ converging to $|\Om|$. Assume that
$\Om_h=\Om_h^1\cup\Om_h^2$ with $\Om_h^i$
open and $\Om_h^1\cap\Om_h^2=\emptyset$. Assume also that $(\Om_h^i)$
converges to an open set
$\Om^i$, $i=1,2$,  in the Hausdorff complementary topology. Then
\begin{itemize}
\item[(i)] $\Om^1\cap\Om^2=\emptyset$,
\item[(ii)] $\Om^1\cup\Om^2=\Om$,
\item[(iii)] $|\Om^i|=\lim_h |\Om^i_h|$, $i=1,2$.
\end{itemize}
In particular, if $\Om_h^0$ is union of connected components of
$\Om_h$ and converges to an open set $\Om^0$
in the Hausdorff complementary topology, then $\Om^0$ is union of
connected components
of $\Om$ and $|\Om^0_h|$ converges to~$|\Om^0|$.
\end{lemma}
\begin{proof} (i) and (ii) are easy
consequences of the convergence in the Hausdorff
complementary topology, while (iii)
follows from Lemma
\ref{convcharac}.
As $\Om_h^0$ is union of connected components of $\Om_h$, then the set
$\Om_h':=\Om_h\setminus\Om_h^0$ is open in the relative topology of
$\Om_h$. Up to a subsequence,
$ \Om_h'$ converges to an open set $\Om'$ in the Hausdorff
complementary topology.
{}From (i) and (ii)  we have that $\Om^0\cap\Om'=\emptyset$, and
$\Om^0\cup\Om'=\Om$;
hence $\Om^0$ is union of connected components of $\Om$. The last
assertion follows from~(iii).
\end{proof}

The following lemma, proved by Bucur and Varchon in \cite{bv1}, will also be
used in the proof of \mbox{Theorem \ref{apcostanti}.}
\begin{lemma}\label{squaresep}
Let $\Om$ be a bounded open set in $\R^2$ and let $a$ and $b$ be two
points in two different connected
components $\Om_a$ and $\Om_b$ of $\Om$, whose distance from $\Om^c$
is greater than $10\delta$ for some $\delta>0$. Let  $U$ be an
    open subset of $\R^2$ such that $U^c$ is connected and
$d_H(U^c,\Om^c)<\delta$.
    Then
there exists $x\in\Om^c$ such that the closed square
$Q(x,9\delta)$, with centre $x$ and side  length $9\delta$,
intersects any curve contained in $U$
and joining the points $a$ and $b$.
\end{lemma}

\begin{proof}[Proof of Theorem \ref{apcostanti}.]
By a density argument, it is sufficient to prove that for every
connected component $\Om^0$ of $\Om$ there exists a
sequence $u_h\in W^{1,p}(\Om_h)$ such that $u_h1_{\Om_h}$ converges
strongly to $ 1_{\Om^0}$ in
$L^p(\R^2)$ and $\nabla u_h1_{\Om_h}$ converges strongly to $ 0$ in
$L^p(\R^2, \R^2)$.
Let $a^0\in \Om^0$
and let $\Om^0_h$ be the connected
component of $\Om_h$ which contains $a^0$ (which is defined for $h$
large enough).
Up to a subsequence, $\Om^0_h$ converges in the Hausdorff complementary
topology to some open set $E$. {}From Lemma \ref{compopalla2}
it follows that
$$E=\bigcup_{i=0}^m\Om^i,$$
where $0\leq m\leq\infty$ and $(\Om^i)$ is a family of connected
components of $\Om$ (including $\Om^0$), and
\begin{equation}\label{stabimise}
|\Om^0_h|\to |E|.
\end{equation}
Let $0<\e<|\Om^0|$ be fixed. There exist a finite integer $n_\e\geq
1$ and an open set $\Om_\e$ such that
\begin{equation}\label{omegaeps}
E=\bigcup_{i=0}^{n_\e}\Om^i\cup\Om_\e,
\end{equation}
where $|\Om_\e|\leq\e$
and $\Om^i\cap\Om_\e=\emptyset$ for every $i\leq n_\e$.

We fix  now a point $a^i$ in each set $\Om^i$.  Let
$\delta>0$ be such that ${\rm dist}(a^i,\Om^c)>10\delta$
for every $i\leq n_\e$.
{}From Lemma \ref{squaresep}, for $h$ big enough
there exist some points $(x^{\delta,i}_h)_{i=1}^{n_\e}$, uniformly
bounded in $\Om^c$,  such that for every $i\leq n_\e$
the square $Q(x^{\delta,i}_h,9\delta)$
   intersects any curve contained in $\Om^0_h$ and joining the points
   $a^0$ and $a^i$. Up
to a subsequence,  we have that $x^{\delta,i}_h\to x^{\delta,i}$ as
$h\to\infty$, for some $x^{\delta,i}\in\Om^c$.
Once again up to a subsequence,  we have that $x^{\delta,i}\to x^{i}$
as $\delta\to 0$, for some $x^{i}\in\Om^c$.
Let
$$K^{\delta,\e}:=\bigcup_{i=1}^{n_\e}\overline{B(x^i,10\delta)}$$
    and, for $i=0,\ldots,n_\e$, let $\Om_h^{\delta,\e,i}$ be the
connected component of $\Om^0_h\setminus K^{\delta,\e}$
    containing $a^i$. As $K^{\delta,\e}\supset  Q(x^{\delta,i}_h,9\delta)$
    for $\delta$ small enough
and $h$ large enough, we have
\begin{equation}\label{ut1}
\Om_h^{\delta,\e,0}\neq\Om_h^{\delta,\e,i}\quad\mbox{ for }i\neq 0.
\end{equation}
    Let $\var^{\delta,\e}$ be the
$C_p$-capacitary potential of $K^{\delta,\e}$, i.e., the solution of
the minimum problem
$$\min\graffe{\int_{\R^2}[|\nabla\var|^p +|\var|^p]\dx :\var\in W^{1,p}(\R^2)
\mbox{, }\var=1\mbox{ }C_p\mbox{-a.e.\ on } K^{\delta,\e}}.$$
We set
\begin{equation}\label{approssimante}
    u_h^{\e,\delta}:=\begin{cases}
1 & \mbox{ in } \Om_h^{\delta,\e,0},\\
\salt
\var^{\delta,\e}  & \mbox{ in }
\Om_h\setminus \Om_h^{\delta,\e,0}.
\end{cases}
\end{equation}
As $\Om_h\cap\partial\Om_h^{\delta,\e,0}\subset K^{\delta,\e}$, we
have that $u_h^{\e,\delta}\in W^{1,p}(\Om_h)$.
    We observe that
$$\|u_h^{\delta,\e}1_{\Om_h}-1_{\Om_h^{\delta,\e,0}}\|_{L^p(\R^2)} +
\|\nabla u_h^{\delta,\e}1_{\Om_h}\|_{L^p(\R^2,\R^2)}
\,\leq\, 2C_p(K^{\delta,\e})^{\frac{1}{p}}.$$
As $C_p(K^{\delta,\e})\leq n_\e C_p(B(0,10\delta))$ and $p\leq 2$, we
conclude that
\begin{equation}\label{approx0}
\limsup_{\delta\to
0}\,\limsup_{h\to\infty}\bigl[\|u_h^{\delta,\e,0}1_{\Om_h}-1_{\Om_h^{\delta,\e}}\|_{L^p(\R^2)} 
+ \|\nabla u_h^{\delta,\e}1_{\Om_h}\|_{L^p(\R^2,\R^2)}\bigr]\,=\,0.
\end{equation}
We claim that
\begin{equation}\label{approxi1}
\limsup_{\delta\to
0}\,\limsup_{h\to\infty}\|u_h^{\e,\delta,0}1_{\Om_h}-1_{\Om^0}\|^p_{L^p(\R^2)}\leq\e,
\end{equation}
and
\begin{equation}\label{approxi2}
\limsup_{\delta\to 0}\,\limsup_{h\to\infty}\|\nabla
u_h^{\e,\delta}1_{\Om_h}\|_{L^p(\R^2,\R^2)}=0,
    \end{equation}
from which the proof of the theorem is achieved by the arbitrariness of $\e$.

It is easy to see that (\ref{approxi2}) follows from (\ref{approx0}),
while (\ref{approxi1}) is
a consequence of (\ref{approx0}) and of the following inequality
\begin{equation}\label{apmeas}
\limsup_{\delta\to
0}\,\limsup_{h\to\infty}|\Om_h^{\delta,\e,0}\triangle\Om^0|\leq\e.
\end{equation}

So, let us prove (\ref{apmeas}).
For every $i=0,\ldots n_\e$, up to a subsequence,
$\Om_h^{\delta,\e,i}$ converges in the Hausdorff
complementary topology, when $h\to\infty$, to some open set
$\Om^{\delta,\e,i}\subset E$.  We observe that $\Om_h^0\setminus
K^{\delta,\e}$ converges to $E\setminus K^{\delta,\e}$
in the Hausdorff complementary topology when $h\to\infty$.
Let $E^{\delta,\e,i}$ be the connected
component  of $E\setminus K^{\delta,\e}$ which contains $a^i$. It is
easy to see that
\begin{equation}\label{ut}
E^{\delta,\e,i}\subset \Om^{\delta,\e,i}.
\end{equation}

Note that, as $\delta\searrow 0$,
$K^{\delta,\e}$ converges decreasingly to the set $\{x_1,\ldots,
x_{n_\e}\}$, $E^{\delta,\e,i}$ converges increasingly
to $\Om^i$, and $\Om^{\delta,\e,i}$ converges increasingly to some
open set $\Om^{\e,i}\subset E$. {}From (\ref{ut}), it follows
that $\Om^i\subset \Om^{\e,i}$.  {}From (\ref{stabimise})
and from Lemma \ref{convcharac} it follows that
\begin{equation}\label{ut3}
|\Om_h^0\setminus K^{\delta,\e}|\to |E\setminus K^{\delta,\e}|.
\end{equation}

By Lemma \ref{compopalla2} applied to $\Om_h^1:=\Om_h^{\delta,\e,0}$ and
$\Om^2_h:=(\Om^0_h\setminus K^{\delta,\e})\setminus \Om_h^{\delta,\e,0}$,
we have that $\Om^{\delta,\e,0}\cap\Om^{\delta,\e,i}=\emptyset$ for
every $i\neq 0$, from which it follows that
    $\Om^{\e,0}\cap\Om^{\e,i}=\emptyset$  and hence
$\Om^{\e,0}\cap\Om^{i}=\emptyset$ for every $1\leq i\leq n_\e$.
    Therefore, there exists an open set $\Om_\e'$, contained in the set $
    \Om_\e$ introduced in (\ref{omegaeps}), such that
$$\Om^{\e,0}=\Om^0\cup\Om_\e'.$$
{}From (\ref{ut3}) and from Lemmas \ref{compopalla2} and
\ref{convcharac}, it follows that
\begin{equation}\label{ut4}
|\Om_h^{\delta,\e,0} \triangle\Om^{\delta,\e,0}|\to 0.
\end{equation}
As $\Om^{\delta,\e,0}\subset\Om^{\e,0}=\Om^0\cup\Om_\e'$, it follows that
\begin{equation}\label{ds1}
|\Om^{\delta,\e,0}_h\setminus\Om^0| \,\leq\,
|\Om^{\delta,\e,0}_h\setminus\Om^{\delta,\e,0}|
+|\Om^{\delta,\e,0}\setminus\Om^0|\,\leq
|\Om^{\delta,\e,0}_h\setminus\Om^{\delta,\e,0}| +
\,|\Om_\e'|,
\end{equation}
and
\begin{equation}\label{ds2}
|\Om^0\setminus\Om^{\delta,\e,0}_h|\,\leq\,|\Om^0
\setminus\Om^{\delta,\e,0}|
+|\Om^{\delta,\e,0}\setminus\Om^{\delta,\e,0}_h|\leq\,|\Om^{\e,0}
\setminus\Om^{\delta,\e,0}|
+|\Om^{\delta,\e,0}\setminus\Om^{\delta,\e,0}_h|.
\end{equation}
As $\Om^{\delta,\e,0}$ converges increasingly to $\Om^{\e,0}$ as
$\delta\to 0^+$,
we have $|\Om^{\e,0}\setminus\Om^{\delta,\e,0}|\to 0$ as $\delta\to 0$. By
(\ref{ut4}),  passing to the limit  in (\ref{ds1}) and (\ref{ds2})
first as $h\to\infty$ and then as $\delta\to
0^+$ we obtain
$$
\limsup_{\delta\to
0}\,\limsup_{h\to\infty}|\Om^{\delta,\e,0}_h\setminus\Om^0|\leq |\Om_\e|
+\limsup_{h\to\infty}|\Om^{\delta,\e,0}_h\setminus\Om^{\delta,\e,0}|\leq\e$$
and
$$\limsup_{\delta\to
0}\,\limsup_{h\to\infty}|\Om^0\setminus\Om^{\delta,\e,0}_h|=0,$$
which give inequality (\ref{apmeas}).
\end{proof}
\begin{remark}\label{Sobpgreat2}
{\rm In the case $p>2$ the stability result for problem 
(\ref{eqgen})  is not true under our hypotheses. Indeed, 
let us  consider
$$S:=[-1,1]\times\{0\},\quad
S_h:=([-1,-1/h]\cup [1/h,1])\times\{0\},$$
$$\Om:= (-1,1)^2\setminus S,
\quad\mbox{ and }\quad\Om_h:= (-1,1)^2\setminus S_h.$$
We set
$$
a(x,\xi):=|\xi|^{p-2}\xi\quad\mbox{ and }\quad 
b(x,\eta):=|\eta|^{p-2}\eta -x_2,
$$
where $x=(x_1,x_2)$. 
Let $u_h$ and $u$ be solutions of problems (\ref{eqgen}) in $\Om_h$ 
and $\Om$ respectively. By the symmetry of $\Om$, the solution $u$ 
 will depend only on $x_2$. Therefore, for every $x=(x_1,x_2)\in\Om$ 
 $$u(x)=\begin{cases}
\mbox{  }w(x_2) &\mbox{ if } x_2\in (0,1),\\
-w(-x_2) &\mbox{ if } x_2\in (-1,0)
\end{cases}$$ 
where $w$ is the solution  of the 
one-dimensional  problem 
\begin{equation}\label{1D}
\begin{cases}
-(|w'|^{p-2}w')'+|w|^{p-2}w=t & \mbox{ in }\, (0,1),\\
\salt
\quad w'(0)= w'(1)=0,& \mbox{  }
\end{cases}
\end{equation}
 which turns out to be of class $C^1([0,1])$.
For every $v\in W^{1,p}(\Om)$, let
  $v^+$ and $v^-$ be the upper and lower traces of $v$ on $S$,
defined as in (\ref{traces}).

If the open set $\Om$
  were stable for problem (\ref{eqgen}) along the sequence $(\Om_h)$,
then $u_h$ would converge
  strongly to $u$ in $W^{1,p}(\Om)$. Hence we would have
that $u_h^+\to u^+$ and $u_h^-\to u^-$ uniformly on $S$
(recall that $p>2$ here). Since $u_h^+(0,0)=u_h^-(0,0)$ by the continuity of 
 $u_h$, we would 
obtain $u^+(0,0)=u^-(0,0)$, which implies that $w(0)=0$. 
Let us prove that this is false. Indeed, by the maximum 
principle we have that $w(t)\geq 0$ 
for every $t\in [0,1]$. Since $w'(0)=0$ and $p>2$, 
we have that $w^{p-1}(t)-t<0$ 
 in a small neighborhood $I$ of $0$ in $[0,1]$. 
So, from equation (\ref{1D}) 
 the function $|w'|^{p-2}w'$ is decreasing in $I$ and hence
 $w'(t)<0$ for every $t\in I$. If $w(0)$ were equal to $0$,
  we would obtain $w(t)<0$ for every $t\in I$, which contradicts 
the fact that  $w(t)\geq 0$ 
for every $t\in [0,1]$. This proves that $w(0)>0$, and hence $\Om$ is not 
 stable for problem (\ref{eqgen}) along the sequence $(\Om_h)$ for $p>2$.

}
\end{remark}
\section{The case of unbounded domains}
We now extend the results of the previous sections to the case of
unbounded domains.
\begin{theorem}\label{msnlim}
Let $(\Om_h)$ be a sequence of open subsets of $\R^2$ such that
$(\Om^c_h)$ converges
to $\Om^c$ in the sense of Kuratowski for some open subset $\Om$.
Assume that, for every $R>0$, $|\Om_h\cap B(0,R)|$ converges to
$|\Om\cap B(0,R)|$
and that the number of connected components of $(\Om_h\cap B(0,R))^c$ is
uniformly bounded with
respect to $h$. Then the sequence of subspaces $X_{\Om_h}$ {\rm
(}resp. $Y_{\Om_h}${\rm )}   converges
    to $X_{\Om}$ {\rm (}resp.\ $Y_{\Om}${\rm )}.
\end{theorem}

\begin{proof}
We  prove only the Mosco convergence of $X_{\Om_h}$ to $X_{\Om}$,
since the convergence of $Y_{\Om_h}$ to
$Y_\Om$ can be proved in the same way.
First of all note that, from the convergence of $\Om^c_h$ to $\Om^c$
in the sense of Kuratowski, it
follows that the sequence $\Om_h\cap B(0,R)$ converges to
$\Om\cap B(0,R)$ in the Hausdorff complementary topology. Moreover,
by the assumptions of the theorem,
we can apply Theorem \ref{moscosobolev} to the
sequence $\Om_h\cap B(0,R)$. So, we get
    that $X_{\Om_h\cap B(0,R)}$ converges to $X_{\Om\cap B(0,R)}$ in the
sense of Mosco.

Now let us prove  ($M_2$) for $X_{\Om_h}$ and $X_\Om$.
Let $(h_k)$ be a sequence of indices
converging to
$\infty$, $(u_k)$ be a sequence, with
$u_k\in W^{1,p}(\Om_{h_k})$ for every $k$, such that $u_k1_{\Om_{h_k}}$
converges weakly
in $L^p(\R^2)$ to a function $\phi$, while $\nabla u_k1_{\Om_{h_k}}$
converges weakly in
$L^p(\R^2,\R^2)$ to a function $\psi$. It follows that
$u_k1_{\Om_{h_k}\cap B(0,R)}$
converges
    to $\phi 1_{B(0,R)}$ weakly in $L^p(\R^2)$, while $\nabla
u_k1_{\Om_{h_k}\cap B(0,R)}$
converges   to $\psi 1_{B(0,R)}$ weakly in $L^p(\R^2,\R^2)$. So, by
property ($M_2$) relative
to the Mosco convergence of $X_{\Om_h\cap B(0,R)}$ to $X_{\Om\cap
B(0,R)}$, there
exists a function $u_R\in W^{1,p}(\Om\cap B(0,R))$ such that $\phi
1_{B(0,R)}=u_R1_{\Om\cap B(0,R)}$
and  $\psi 1_{B(0,R)}=\nabla u_R1_{\Om\cap B(0,R)}$ a.e.\ in $\R^2$.
Since $R$ is arbitrary, it is easy to construct $u\in W^{1,p}(\Om)$
such that $\phi =u1_\Om$ and $\psi=\nabla u1_\Om$ a.e.\ in $\R^2$.

Let us prove property ($M_1$). Let $u\in W^{1,p}(\Om)$ and let
$\e>0$. There exists $R_\e>0$ such that
$$\int_{\Om\setminus B(0,R_\e)}[|u|^p+|\nabla u|^p]\dx\leq\e.$$ By
property ($M_1$) relative
to the Mosco convergence of $X_{\Om_h\cap B(0,R_\e+1)}$ to
$X_{\Om\cap B(0,R_\e+1)}$, there exists a sequence
$w^\e_h\in W^{1,p}(\Om_h\cap B(0,R_\e+1))$ such that
$w^\e_h1_{\Om_h\cap B(0,R_\e+1)}$ converges strongly
to $u 1_{\Om\cap B(0,R_\e+1)}$ in $ L^p(\R^2)$ and
$\nabla w^\e_h1_{\Om_h\cap B(0,R_\e+1)}$
converges strongly to $\nabla u 1_{\Om\cap B(0,R_\e+1)}$  in $L^p(\R^2,\R^2)$.
    Let $\var_\e\in C^1_c(B(0,R_\e+1))$ such that
    $0\leq\var_\e\leq 1$,\,\,$\var_\e=1$ in $B(0,R_\e)$, and
$\|\nabla\var_\e\|_\infty\leq C$.
Now we set $u^\e_h:=\var _\e w^\e_h$. By construction
$u_h^\e\in W^{1,p}(\Om_h)$, $u_h^\e1_{\Om_h}\to\var_\e u1_\Om$ strongly
in $L^p(\R^2)$, and $\nabla u_h^\e1_{\Om_h}\to\var_\e\nabla u1_\Om
+u\nabla\var_\e1_\Om$ strongly in $L^p(\R^2,\R^2)$. On the other hand
$$
\limsup_{h\to\infty}\int_{\R^2}|u^\e_h1_{\Om_h}-u1_{\Om}|^p
+ |\nabla u^\e_h1_{\Om_h}-\nabla u1_{\Om}|^p\dx\,\leq \,2^{p-1}(C^p+1)\e,
$$
$\var_\e u1_\Om\to u1_\Om$ strongly
in $L^p(\R^2)$, and $\var_\e\nabla u1_\Om
+u\nabla\var_\e1_\Om\to \nabla u1_\Om$ strongly in $L^p(\R^2,\R^2)$.
   Therefore, we can construct a sequence $u_h\in
W^{1,p}(\Om_h)$ which satisfies ($M_1$) by a standar argument on
double sequences.
\end{proof}
\section{Problems  with Dirichlet boundary
conditions}
    In this section we study the Mosco convergence of Sobolev
and Deny-Lions
spaces with prescribed Dirichlet conditions on part of the boundary.

    Let $A\subset\R^2$ be a bounded open set with Lipschitz
   boundary $\partial A$, and let $\partial_DA$ be a
   relatively open subset  of $\partial A$ with a finite number
of connected components.  For every compact set $K\subset\overline A$,
    for every $g\in W^{1,p}(A)$, and for every pair of function $a$
    and $b$ satisfying the properties (\ref{a1})--(\ref{a3}), we consider
    the solutions $u$ and $v$ of the mixed problems
\begin{equation}\label{eqmix1} \begin{cases} -\,{\rm
div}\,a(x,\nabla u)+  b(x,u)=0 & \mbox{ in }\quad A\setminus K,\\
    \salt
\mbox{ }u=g & \mbox{ on
}\quad \partial_DA\setminus K,\\
\salt
\mbox{ }a(x,\nabla u)\cdot\nu=0 & \mbox{ on
}\quad \partial (A\setminus K)\setminus (\partial_DA\setminus K),
\end{cases}
\end{equation}
and
\begin{equation}\label{eqmix2} \begin{cases} -\,{\rm
div}\,a(x,\nabla v)=0 & \mbox{ in }\quad A\setminus K,\\
    \salt
\mbox{ }v=g & \mbox{ on
}\quad \partial_DA\setminus K,\\
\salt
\mbox{ }a(x,\nabla v)\cdot\nu=0 & \mbox{ on
}\quad \partial (A\setminus K)\setminus (\partial_DA\setminus K).
\end{cases}
\end{equation}
Let $(K_h)$ be a sequence of compact subsets of $\overline A$, let
$(g_h)$ be a sequence in $W^{1,p}(A)$, and let $(u_h)$ be the sequence
   of the solutions of problems (\ref{eqmix1})  corresponding to
   $K_h$ and $g_h$.
\begin{definition}\label{stabledc}
We say that the pair $(K,g)$ is stable for the mixed  problems
(\ref{eqmix1}) along the sequence $(K_h,g_h)$ if for
every pair of functions $a$, $b$
satisfying (\ref{a1})--(\ref{a3}) the sequence
$(u_h1_{K^c_h})$ converges to $u1_{K^c}$ strongly in $L^p(A)$
and the sequence $(\nabla u_h1_{K^c_h})$ converges to $\nabla u1_{K^c}$
strongly in
$L^p(A,\R^2)$.

The stability for problems (\ref{eqmix2}) is defined in a similar way
by using only the convergence of the gradients {\rm (}as in Definition
~\ref{stableaux}{\rm )}.
\end{definition}

The stabilty for problems (\ref{eqmix2}) has been recently studied
in \cite{dmro} in the case $a(x,\xi)=\xi$. In this section we will
study the stability in the general case by using again the notion of Mosco
convergence.

We set
$$W^{1,p}_g(A\setminus K,\partial_DA\setminus K):=\{u\in W^{1,p}(A\setminus K):
u=g\mbox{ on
} \partial_DA\setminus K\},$$
and
$$
L^{1,p}_g(A\setminus K,\partial_DA\setminus K):=\{u\in L^{1,p}(A\setminus K):
u=g\mbox{ on
} \partial_DA\setminus K\},$$
where the equality $u=g$ on $\partial_DA\setminus K$ is intended in
the usual sense
of traces.

As in Section \ref{snp} the space $W^{1,p}_g(A\setminus
K,\partial_DA\setminus K)$
will be identified with the closed linear subspace $X_K^g(A)$ of $L^p(A)\times
L^p(A,\R^2)$ defined by
\begin{equation}\label{XAg}
X_K^g(A):=\{(u1_{K^c},\nabla u1_{K^c}): u\in W^{1,p}_g(A\setminus
K,\partial_DA\setminus K)\}.
\end{equation}

    For problem (\ref{eqmix2}), we consider in $L^{1,p}_g(A\setminus
K,\partial_DA\setminus K)$ the
equivalence relation $\sim$ defined in (\ref{defreleq}).
Note that in this case $v_1\sim v_2$ if and only if $v_1=v_2$ a.e.\
in those connected components
of $A\setminus K$ whose boundary intersects $\partial_DA\setminus K$
and $\nabla v_1=\nabla v_2$ a.e.\ in the other connected components
of $A\setminus K$.
The corresponding quotient space,
denoted by $L^{1,p}_g(A\setminus K,\partial_DA\setminus K)/_\sim$,
will be identified with the closed linear subspace
$Y_K^g(A)$ of $ L^p(A,\R^2)$ defined by
\begin{equation}\label{YAg}
Y_K^g(A):=\{\nabla u1_{K^c}: u\in L^{1,p}_g(A\setminus
K,\partial_DA\setminus K)\}.
\end{equation}

Let $K_h$, $K$ be compact subsets of $\overline A$ and let $g_h$,
$g\in W^{1,p}(A)$.
Let $X_{K_h}^{g_h}(A)$ and $X_K^g(A)$ be the corresponding subspaces
defined by (\ref{XAg}).
We recall that $X_{K_h}^{g_h}(A)$ converges to $X_K^g(A)$ in the
sense of Mosco if the following two properties hold:
\begin{itemize}
\item[($M_1''$)] for every $u\in W^{1,p}_g(A\setminus
K,\partial_DA\setminus K)$,
there
exists a sequence $u_h\in W^{1,p}_{g_h}(A\setminus
K_h,\partial_DA\setminus K_h)$
such that
$u_h1_{K^c_h}$ converges strongly to $u1_{K^c}$ in  $L^p(A)$ and
$\nabla u_h1_{K^c_h}$
converges
strongly to $\nabla u1_{K^c}$ in  $L^p(A,\R^2)$;
\item[($M_2''$)] if $(h_k)$ is a sequence of indices converging to
$\infty$, $(u_k)$ is a sequence, with
$u_k\in W^{1,p}_{g_{h_k}}(A\setminus K_{h_k}, \partial_DA\setminus
K_{h_k})$ for every
$k$, such that $u_k1_{K^c_{h_k}}$ converges
weakly in $L^p(A)$ to a function $\phi$, while $\nabla u_k1_{K^c_{h_k}}$
converges weakly in  $L^p(A,\R^2)$ to a function
$\psi$,
   then there exists $u\in W^{1,p}_g(A\setminus
K,\partial_DA\setminus K)$ such that $\phi=u1_{K^c}$ and
$\psi=\nabla u1_{K^c}$ a.e.\ in $A$.
\end{itemize}
Analogously, the convergence of $Y_{K_h}^{g_h}(A)$ to $Y_K^g(A)$
in the sense of Mosco can be characterized by using only the
convergence of the extensions of the
gradients.
\begin{remark}\label{rem}
{\rm As in Section \ref{MC} we can prove that the Mosco convergence of
$X_{K_h}^{g_h}(A)$ to $X_K^g(A)$  (resp.\ of $Y_{K_h}^{g_h}(A)$ to
$Y_K^g(A)$) is equivalent to the stability of
$(K,g)$ for the mixed  problems
(\ref{eqmix1})  (resp.\ (\ref{eqmix2}) along the sequence $(K_h,g_h)$).}
\end{remark}

The following theorem is the main result of this section.

\begin{theorem}\label{moscodati}
Let $A$ be a bounded open subset of $\R^2$ with Lipschitz
boundary $\partial A$ and let $\partial_DA$ be a relatively open
subset of $\partial A$ with a finite
number of connected components.  Let $(g_h)$ be a sequence in
$W^{1,p}(A)$ converging strongly to a function $g$
    in $W^{1,p}(A)$,  and let $(K_h)$ be a sequence
    of compact subsets of $\overline A$
converging to a set
    $K$ in the Hausdorff metric. Assume that $|K_h |$ converges to $|K|$
and that the sets $K_h$ have a uniformly bounded
number of connected components.
Then $X_{K_h}^{g_h}(A)$ converges to $X_K^g(A)$ 
{\rm (}resp.\ $Y_{K_h}^{g_h}(A)$ converges to $Y_K^g(A)${\rm )}
    in the sense of Mosco.
    \end{theorem}

\begin{proof}  The main idea of this proof is due to Chambolle.
      Let us first prove the Mosco convergence of $X_{K_h}^{g_h}(A)$
    to $X_K^g(A)$. Let $\Sigma$ be an open ball in $\R^2$ such that
$\overline A\subset\Sigma$. Let $\tilde g_h$, $\tilde g\,\in
W^{1,p}(\Sigma)$
be  extensions of $g_h$ and $g$ to $\Sigma$ such that $\tilde g_h$
converges to $\tilde g$ strongly in $W^{1,p}(\Sigma)$.
We set
$$\Om_h:=\Sigma\setminus(K_h\cup (\partial
A\setminus\partial_DA))\quad\mbox{ and }
\quad\Om:=\Sigma\setminus(K\cup (\partial A\setminus\partial_DA)).$$
Note that $\Om_h$ and $\Om$ satisfy the assumptions of Theorem
\ref{moscosobolev}.
Let us prove property $(M2'')$. Let $(h_k)$ be a sequence of indices that
tends to $\infty$, and $u_k\in W^{1,p}_{g_{h_k}}(A\setminus
K_{h_k},\partial_DA\setminus K_{h_k})$ such that $u_k1_{K^c_{h_k}}$
converges weakly to $\phi$ in $L^p(A)$ and
$\nabla u_k1_{K^c_{h_k}}$ converges weakly to $\psi$ in $L^p(A,\R^2)$.
Let $\tilde u_k$ be the extension of $u_k$ defined by
$$
\tilde u_k:=\begin{cases}
u_k1_{K^c_{h_k}} &\mbox{ in }\, A,\\
\tilde g_{h_k} &\mbox{ in } \Sigma\setminus A,\end{cases}
$$ and let $\tilde\phi$ and $\tilde\psi$ be defined by
$$
\tilde \phi:=\begin{cases}
\phi &\mbox{ in }\, A,\\
\tilde g &\mbox{ in } \Sigma\setminus A,\end{cases}
\qquad\qquad\mbox{ and }\qquad\qquad
\tilde \psi:=\begin{cases}
\psi &\mbox{ in }\, A,\\
\nabla\tilde g &\mbox{ in } \Sigma\setminus A.\end{cases}$$
    As $u_k=g_{h_k}$ on $\partial_DA\setminus K_{h_k}$,
we have $\tilde u_k\in W^{1,p}(\Om_{h_k})$.
Since $\tilde u_k1_{\Om_{h_k}}$ converges to $\tilde\phi 1_\Sigma$
weakly in $L^p(\R^2)$ and $\nabla \tilde u_k1_{\Om_{h_k}}$
converges to $\tilde\psi 1_\Sigma$ weakly in $L^p(\R^2,\R^2)$,
    by Theorem \ref{moscosobolev}
we conclude that there exists $\tilde u\in W^{1,p}(\Om)$ such that
$\tilde\phi 1_\Sigma =\tilde u 1_\Om$ and $\tilde\psi 1_\Sigma
=\nabla\tilde u 1_\Om$.
Let $u$ be the restriction of $\tilde u$ to $A\setminus K$.
Then $u\in W^{1,p}(A\setminus K)$ and we have that
$\phi=u1_{K^c}$ and $\psi=\nabla u1_{K^c}$ a.e.\ in $A$.
As $\tilde u\in W^{1,p}(\Om)$, the traces of
$\tilde u$ on both sides of $\partial_DA\setminus K$ coincide.
Since $\tilde u=\tilde g$ a.e.\ in $\Sigma\setminus A$, we conclude
that $u=\tilde g=g$ in the sense of traces on $\partial_DA\setminus K$.
Therefore $u\in W^{1,p}_g(A\setminus K,\partial_DA\setminus K)$.

Now we prove property ($M_1''$).
Let $u\in W^{1,p}_g(A\setminus K,\partial_DA\setminus K)$.
Let $\tilde u$ be the extension of $u$ defined by
\[ \tilde u=
\begin{cases}
u1_{K^c} & \text{ a.e.\ in } A,\\
\tilde g & \text{ a.e.\ in } \Sigma\setminus A.
\end{cases}
\]
As $u=g$ on $\partial_DA\setminus K$, we have that $\tilde u\in
W^{1,p}(\Om)$. By Theorem \ref{moscosobolev} there exists
a sequence $\tilde u_h\in  W^{1,p}(\Om_h)$ such
that $\tilde u_h1_{\Om_h}$ converges to $\tilde u1_\Om$  strongly in
$L^p(\R^2)$
and   $\nabla \tilde u_h1_{\Om_h}$ converges to $\nabla\tilde u1_\Om
$  strongly
in $L^p(\R^2,\R^2).$
We consider the function
$$\phi_h:=( \tilde u_h - \tilde g_h)|_{\Sigma\setminus\overline A}.$$
By construction, $\phi_h\to 0$ strongly in
$W^{1,p}(\Sigma\setminus\overline A)$. Therefore there exists a sequence
$v_h\in W^{1,p}(\Sigma)$, converging to $0$ strongly in
$W^{1,p}(\Sigma)$,
such that $v_h|_{\Sigma\setminus\overline A}=\phi_h$ a.e.\ in
$\Sigma\setminus\overline A$.
    We set
$$u_h: = \bigl( \tilde u_h  - v_h \bigr)|_{A\setminus K_h}.$$
By construction, we have that $u_h \in
W^{1,p}(A \setminus K_h)$ and $u_h = g_h$ in the sense of traces on
$\partial_D A\setminus K_h$.
Moreover, we have that
$ u_h1_{K^c_h}$ converges  to $u1_{K^c}$ strongly in $L^p(A)$ and
$\nabla u_h1_{K^c_h}$
converges to $\nabla u1_{K^c}$ strongly in $L^p(A,\R^2)$.
\vskip .2truecm
    Now let us prove that $Y_{K_h}^{g_h}(A)$ converges to $Y_K^g(A)$
in the sense of Mosco. Property ($M_2''$) is obtained arguing as in
   \cite[Lemma 4.1]{dmro}. So, let us prove ($M_1''$).
Let $u\in L^{1,p}_g(A\setminus K,\partial_DA\setminus K)$.
We set for every $k\in\N$
$$u^k:=(g-k)\vee u\wedge(g+k).$$
Then $ u^k\in W^{1,p}_g(A\setminus
K,\partial_DA\setminus K)$ and $\nabla u^k1_{K^c}\to\nabla u1_{K^c}$
   strongly in $L^p(A,\R^2)$.
{}From property ($M_1''$) proved above for the Mosco convergence
    of $X_{K_h}^{g_h}(A)$ to $X_K^g(A)$, for every $k$ there exists
$u^k_h\in W^{1,p}_{g_h}(A\setminus
K_h,\partial_DA\setminus K_h)$ such that
$\nabla u^k_h1_{K_h^c}\to\nabla u^k1_{K^c}$
in $L^p(A,\R^2)$.
Hence, by a standard argument on double sequences, we obtain a
sequence of indices $k_h$
converging to $\infty$ such that, setting $u_h:=u^{k_h}_h$, we get
$\nabla u_h1_{K_h^c}\to\nabla u1_{K^c}$ in $L^p(A,\R^2)$.
\end{proof}

\medskip
\centerline{\sc Acknowledgements}
This work is part of the
European Research Training Network ``Homogenization and Multiple
Scales'' under contract HPRN-2000-0010, and of the Research \
Project ``Calculus of Variations''
supported by SISSA and by
the Italian Ministry of Education, University, and Research.


\begin{thebibliography}{99}
{\frenchspacing
\bibitem{BO-Mu-Pu}{Boccardo L., Murat F., Puel J.P.}: {Existence of
bounded solutions for
nonlinear elliptic unilateral problems}, {\it Ann. Mat. Pura Appl. (4)}
    {\bf 152}  (1988), 183-196.
\bibitem{bv1} {Bucur D., Varchon N.}: {Boundary variation for a 
Neumann problem}, {\it Ann. Scuola Norm. Sup. Pisa Cl. Sci. (4)}
 {\bf 29} (2000), 807-821.
\bibitem{bv2} {Bucur D., Varchon N.}: {A duality approach for the
boundary variation of Neumann problems}, Preprint Univ.
Franche-Comt\'e, 2000.
\bibitem{bt} {Bucur D., Trebeschi P.}: {Shape optimization problem
governed by a nonlinear state equation}, {\it Proc. Roy. Soc.
Edinburgh Sect. A} {\bf 128} (1998), 949-963.
\bibitem{buzo1} {Bucur, D., Zolesio, J.-P.}:
{Shape optimization for elliptic problems under
Neumann boundary conditions,}
Bouchitt\'e, G. (ed.) et al., {\it Calculus of Variations,
Homogenization and Continuum Mechanics
(CIRM-Luminy, Marseille, 1993).} Singapore: World Scientific, Ser.
{\it Adv. Math. Appl. Sci.} {\bf 18} (1994), 117-129.
\bibitem{buzo2} {Bucur, D., Zolesio, J.-P.}:
Continuit\'e par rapport au domaine dans le probl\`eme de Neumann,
{\it C. R. Acad. Sci. Paris S\'er. I} {\bf 319} (1994), 57-60.
\bibitem{buzo3} {Bucur, D., Zolesio, J.-P.}:
{Shape continuity for Dirichlet-Neumann problems,}
Chipot, M. (ed.) et al., {\it Progress in Partial Differential Equations:
the Metz Surveys 4. (Metz, 1994-95)}.
Pitman Res. Notes in Math. 345, Longman, Harlow (1996), 53-65.
\bibitem{chdo}{Chambolle A., Doveri F.}: {Continuity of Neumann
linear elliptic problems on varying two-dimensional boun\-ded open
sets}, {\it Comm. Partial Differential Equations} {\bf 22}
(1997), 811-840.
\bibitem{chen} {Chenais D.}: {On the existence of a solution in
domain identification problem}, {\it J. Math Anal. Appl.} {\bf 52}
(1975), 189-289.
\bibitem{con} {Conca C.}: {The Stokes sieve problem,}
{\it Comm. Appl. Num. Meth.} {\bf 4} (1988), 113-121.
\bibitem{cor} {Cortesani G.}: {Asymptotic behaviour of a sequence of
Neumann problems}, {\it Comm. Partial Differential Equations} {\bf
22} (1997), 1691-1729.
\bibitem{dmro} {Dal Maso G., Toader R.}: {A model for the
quasi-static growth of brittle fractures: existence and approximation
results,} {\it Arch. Ration. Mech. Anal.} {\bf 162} (2002), 101-135.
\bibitem{daml} {Damlamian A.}: {Le probl\`eme de la passoire de
Neumann}, {\it Rend. Sem. Mat. Univ. Politec. Torino} {\bf 43} (1985), 427-450.
\bibitem{delv}{Del Vecchio T.}: {The thick Neumann's sieve,}
{\it Ann. Mat. Pura Appl. (4)} {\bf 147} (1987), 363-402.
\bibitem{deli} {Deny J., Lions J.L.}: {Les espaces du type Beppo
Levi,} {\it Ann. Inst. Fourier (Grenoble)} {\bf 5} (1953), 305-370.
\bibitem{ebpo} {Ebobisse F., Ponsiglione M.}:
{Limit of some mixed
type boundary value problems for the $p$-laplacian equation in
nonsmooth domains}, forthcoming.
\bibitem{EG}Evans L.C., Gariepy R.F.: {\it Measure Theory and Fine Properties
of Functions}, CRC Press, Boca Raton, 1992.
\bibitem{falc} {Falconer K.J.}: {\it The Geometry of Fractal Sets},
Cambridge University Press, Cambridge, 1985.
\bibitem{hekima} {Heinonen J., Kipelanen T., Martio O.}: {\it
Nonlinear Potential Theory of Degenerate Elliptic Equations},
Clarendon Press, Oxford, 1993.
\bibitem{lions} {Lions J.L.}: {\it Quelques m\'ethodes de
r\'esolution des probl\`emes aux limites non lin\'eaires}, Dunod
Paris; Gauthier-Villars, Paris, 1968.
\bibitem{maz1} {Maz'ya V.G.}: {\it Sobolev Spaces}, Springer-Verlag,
Berlin, 1985.
\bibitem{mosco} {Mosco U.}: {Convergence of convex sets and of
solutions of variational inequalities}, {\it Adv. Math.} {\bf 3}
(1969), 510-585.
\bibitem{mur} {Murat F.}: {The Neumann sieve}, {\it Nonlinear Variational
Problems (Isola d'Elba, 1983)}, Marino A. {\it et al}., editor,
Res. Notes in Math. 127, Pitman, London, (1985), 24-32.
\bibitem{sver} {\v Sver\'ak V.}: { On optimal shape design}, {\it J.
Math. Pures Appl. (9)} {\bf 72} (1993), 537-551.
\bibitem{ziem} {Ziemer W.} {\it Weakly Differentiable Functions},
Springer-Verlag, Berlin, 1989.
\par
}
\end{thebibliography}
\end{document}